\documentstyle[12pt]{article}
\begin{document}
\pagestyle{empty}
\renewcommand{\thefootnote}{\fnsymbol{footnote}}
\def\lsim{\raise0.3ex\hbox{$<$\kern-0.75em\raise-1.1ex\hbox{$\sim$}}}
\def\gsim{\raise0.3ex\hbox{$>$\kern-0.75em\raise-1.1ex\hbox{$\sim$}}}
\def\noi{\noindent}
\def\sq{\hbox {\rlap{$\sqcap$}$\sqcup$}}
\def\R{ {\rm R \kern -.31cm I \kern .15cm}}
\def\C{ {\rm C \kern -.15cm \vrule width.5pt \kern .12cm}}
\def\Z{ {\rm Z \kern -.27cm \angle \kern .02cm}}
\def\N{ {\rm N \kern -.26cm \vrule width.4pt \kern .10cm}}
\def\1{{\rm 1\mskip-4.5mu l} }
\def\lsim{\raise0.3ex\hbox{$<$\kern-0.75em\raise-1.1ex\hbox{$\sim$}}}
\def\gsim{\raise0.3ex\hbox{$>$\kern-0.75em\raise-1.1ex\hbox{$\sim$}}}
\vbox to 2 truecm {}
\centerline{\Large \bf Long Range Scattering and Modified}
\vskip 3 truemm  
\centerline{\Large \bf Wave Operators for some Hartree Type Equations II\footnote{Work supported in part by NATO Collaborative Research Grant 972231.}}

\vskip 1 truecm
\centerline{\bf J. Ginibre}
\centerline{Laboratoire de Physique Th\'eorique\footnote{Unit\'e Mixte de Recherche (CNRS) UMR
8627.}}  \centerline{Universit\'e de Paris XI, B\^atiment 210,
F-91405 Orsay Cedex, France}

\vskip 5 truemm
\centerline{\bf G. Velo\footnote{Permanent address : Dipartimento di Fisica, Universit\`a di Bologna and INFN, Sezione di
Bologna, Italy. }}
\centerline{Laboratoire d'Analyse Num\'erique et E.D.P.\footnote{Unit\'e Mixte de Recherche (CNRS)
UMR 8628.}}  \centerline{Universit\'e de Paris XI, B\^atiment 425,
F-91405 Orsay Cedex, France}

\vskip 1 truecm 
\begin{abstract}
We study the theory of scattering for a class of Hartree type equations with long range
interactions in space dimension $n \geq 3$, including Hartree equations with potential $V(x) =
\lambda |x|^{- \gamma}$. For $0 < \gamma \leq 1$ we prove the existence of
modified wave operators with no size restriction on the data and we determine the asymptotic
behaviour in time of solutions in the range of the wave operators, thereby extending the results of
a previous paper which covered the range $1/2 < \gamma < 1$.
   \end{abstract}
\vskip 3 truecm
\noi AMS Classification : Primary 35P25. Secondary 35B40, 35Q40, 81U99.  \par
\noi Key words : Long range scattering, modified wave operators, Hartree equation. 
\vskip 1 truecm

\noindent LPT Orsay 99-14 \par
\noindent February 1999 \par

  \newpage
\pagestyle{plain}
\baselineskip=24 pt

\section{Introduction}
\hspace*{\parindent} This is the second paper where we study the theory of scattering and more
precisely the existence of modified wave operators for a class of long range Hartree type equations
$$i \partial_t u  + {1 \over 2} \Delta u = \widetilde{g} (|u|^2)u \eqno(1.1)$$
\noi where $u$ is a complex function defined in space time ${I\hskip-1truemm R}^{n+1}$, $\Delta$
is the Laplacian in ${I\hskip-1truemm R}^n$, and
$$\widetilde{g} (|u|^2) = \lambda t^{\mu - \gamma} \ \omega^{\mu - n} \ |u|^2 \eqno(1.2)$$
\noi with $\omega = (- \Delta)^{1/2}$, $\lambda \in {I\hskip-1truemm R}$, $0 < \gamma \leq 1$ and $0
< \mu < n$. The operator $\omega^{\mu - n}$ can also be represented by the convolution in $x$
$$\omega^{\mu - n} \ f = C_{n, \mu} \ |x|^{- \mu} * f \eqno(1.3)$$
\noi so that (1.2) is a Hartree type interaction with potential $V(x) = C |x|^{-
\mu}$. The more standard Hartree equation corresponds to the case $\gamma = \mu$. In that case,
the nonlinearity $\widetilde{g}(|u|^2)$ becomes
$$\widetilde{g}(|u|^2) = V * |u|^2 = \lambda |x|^{-\gamma} * |u|^2 \eqno(1.4)$$
\noi with a suitable redefinition of $\lambda$. \par

A large amount of work has been devoted to the theory of scattering for the Hartree equation
(1.1) with nonlinearity (1.4) as well as with similar nonlinearities with more general potentials.
As in the case of the linear Schr\"odinger equation, one must
distinguish the short range case, corresponding to $\gamma > 1$, from the long range case
corresponding to $\gamma \leq 1$. In the short range case, it is known that the (ordinary) wave
operators exist in suitable function spaces for $\gamma > 1$ \cite{11r}. Furthermore for repulsive
interactions, namely for $\lambda \geq 0$, it is known that all solutions in suitable spaces admit
asymptotic states in $L^2$ for $\gamma > 1$, and that asymptotic completeness holds for $\gamma >
4/3$ \cite{10r}. In the long range case $\gamma \leq 1$, the ordinary wave operators are known not
to exist in any reasonable sense \cite{10r}, and should be replaced by
modified wave operators including a suitable phase in their definition, as is the case for the
linear Schr\"odinger equation. A well developed theory of long range scattering exists for the
latter. See for instance \cite{1r} for a recent treatment and for an extensive
bibliography. In contrast with that situation, only partial results are available for the
Hartree equation. On the one hand, the existence of
modified wave operators has been proved in the critical case $\gamma = 1$ for small solutions
\cite{2r}. On the other hand, it has been shown, first in the critical case $\gamma = 1$
\cite{6r,9r} and then in the whole range $0 < \gamma \leq 1$ \cite{5r,7r,8r} that the global
solutions of the Hartree equation (1.1) (1.3) with small initial data exhibit an asymptotic
behaviour as $t \to \pm \infty$ of the expected scattering type characterized by scattering states
$u_{\pm}$ and including suitable phase factors that are typical of long range scattering. In
particular, in the framework of scattering theory, the results of \cite{5r,7r,8r} are closely
related to the property of asymptotic completeness for small data. \par

In a previous paper with the same title \cite{4r}, hereafter referred to as I, we proved the
existence of modified wave operators for the equation (1.1) (1.2), and we gave a description of the
asymptotic behaviour in time of solutions in the ranges of those operators, with no size restriction
on the data, in suitable spaces and for $\gamma$ in the range $1/2 < \gamma < 1$. The method is an
extension of the energy method used in \cite{5r,7r,8r}, and uses in particular the equations
introduced in \cite{7r} to study the asymptotic behaviour of small solutions. The spaces of initial
data, namely in the present case of asymptotic states, are Sobolev spaces of finite order similar to
those used in \cite{8r}. The present paper is devoted to the extension of the previous results to
the whole range $0 < \gamma \leq 1$. The methods used here are natural extensions of those used in
I. They require in particular the same restrictions on $\mu$ and $n$, in particular $\mu \leq n -
2$ and $n \geq 3$. We refer to the introduction of I for a discussion of those conditions. \par

The construction of the modified wave operators is too complicated to allow for a more precise
statement of results at this stage, and will be described in Section 2 below, which is a summary
and continuation of Section 2 of I. That construction involves the study of the same auxiliary
system of equations as in I, for a new function $w$ and a phase $\varphi$ instead of the original
function $u$, and relies as a preliminary step on the construction of local wave operators in a
neighborhood of infinity for that system. That step requires the definition of a modified
asymptotic dynamics which is significantly more complicated than that used in I. \par

We now give a brief outline of the contents of this paper. A more detailed description of the
technical parts will be given at the end of Section 2. After collecting some notation and
preliminary estimates in Section 3 and recalling from I some preliminary results on the auxiliary
system in Section 4, we define and study the asymptotic dynamics in Section 5. We then study the
asymptotic behaviour of solutions for the auxiliary system in
Section 6. In particular we essentially construct local wave operators at infinity
for that system. We then come back from the auxiliary system to the original equation
(1.1) for $u$ and construct the wave operators for the latter in Section 7, where the final
result will be stated in Proposition 7.5. \par

We have tried to make this paper as self-contained as possible and at the same time to keep
duplication with I to a minimum. Duplication occurs in the beginning of Section 3 and in Section
4 where we recall estimates and results from I. On the other hand, Sections 6 and 7 follow the
same pattern as Sections 5, 6 and 7 of I, with the appropriate changes needed to handle the more
general situation.\par

We conclude this section with some general notation which will be used freely throughout this
paper. We denote by $\parallel \cdot \parallel_r$ the norm in $L^r \equiv L^r ({I\hskip-1truemm
R}^n)$. For any interval $I$ and
any Banach space $X$, we denote by ${\cal C}(I, X)$ the space of strongly continuous functions
from $I$ to $X$ and by $L^{\infty}(I, X)$ (resp. $L_{loc}^{\infty}(I, X))$ the space of
measurable essentially bounded (resp. locally essentially bounded) functions from $I$ to $X$. For
real numbers $a$ and $b$, we use the notation $a \vee b = {\rm Max} (a, b)$, $a \wedge b = {\rm
Min}(a,b)$ and $[a] =$ integral part of $a$. In the estimates of solutions of the relevant
equations, we shall use the letter $C$ to denote constants, possibly different from an
estimate to the next, depending on various parameters such as $\gamma$, but not on the
solutions themselves or on their initial data. Those constants will be bounded in $\gamma$
for $\gamma$ away from zero. We shall use the notation $A(a_1,a_2,\cdots)$ for estimating
functions, also possibly different from an estimate to the next, depending in addition on
suitable norms $a_1, a_2, \cdots$ of the solutions or of their initial data. Finally Item (p, q) of I
will be referred to as Item (I. p. q). Additional notations will be given at the beginning of
Section 3. \par

In all this paper, we assume that $n \geq 3$, $0 < \mu \leq n - 2$ and $0 < \gamma \leq 1$.

\section{Heuristics}
\hspace*{\parindent} In this section, we discuss in heuristic terms the construction of the modified
wave operators for the equation (1.1), as it will be performed in this paper. That construction is an
extension of that performed in I in the special case $\gamma > 1/2$, and we refer to Section I.2 for
a more detailed introduction and for general background. \par

The problem that we want to address is that of classifying the possible asymptotic behaviours of
the solutions of (1.1) by relating them to a set of model functions ${\cal V} = \{v = v(u_+)\}$
parametrized by some data $u_+$ and with suitably chosen and preferably simple asymptotic
behaviour in time. For each $v \in {\cal V}$, one tries to construct a solution $u$ of (1.1) such
that $u(t)$ behaves as $v(t)$ when $t \to \infty$ in a suitable sense. The map $\Omega : u_+ \to
u$ thereby obtained classifies the asymptotic behaviours of solutions of (1.1) and is a
preliminary version of the wave operator for positive time. A similar question can be asked for
$t \to - \infty$. From now on we restrict our attention to positive time. \par

In the short range case corresponding to $\gamma > 1$ in (1.1), the previous scheme can be
implemented by taking for ${\cal V}$ the set ${\cal V} = \{ v = U(t) u_+\}$ of solutions of the
equation
$$i \partial_t v + {1\over 2} \Delta v = 0 \quad , \eqno(2.1)$$

\noi with $U(t)$ being the unitary group
$$U(t) = \exp \left ( i(t/2)\Delta \right ) \quad . \eqno(2.2)$$

\noi The initial data $u_+$ for $v$ is called the asymptotic state for $u$. \par

In the long range case corresponding to $\gamma \leq 1$ in (1.1) (1.2), the previous set is known to
be inadequate and has to be replaced by a better set of model functions obtained by modifying the
previous ones by a suitable phase. The modification that we use requires additional structure of
$U(t)$. In fact $U(t)$ can be written as
$$U(t) = M(t) \ D(t) \ F \ M(t) \eqno(2.3)$$

\noi where $M(t)$ is the operator of multiplication by the function 
$$M(t) = \exp \left ( i  x^2/2t \right ) \quad , \eqno(2.4)$$

\noi $F$ is the Fourier transform and $D(t)$ is the dilation operator defined by 
$$\left ( D(t) \ f \right ) (x) = (it)^{-n/2} \ f(x/t) \quad . \eqno(2.5)$$

\noi Let now $\varphi^{(0)} = \varphi^{(0)}(x,t)$ be a real function of space time and let
$z^{(0)}(x, t) = \exp (-i \varphi^{(0)}(x,t))$. We replace $v(t) = U(t) u_+$ by the modified free
evolution \cite{12r} \cite{13r}
$$v(t) = M(t) \ D(t) \ z^{(0)}(t) \ w_+   \eqno(2.6)$$

\noi where $w_+ = Fu_+$. In order to allow for easy comparison of $u$ with $v$, it is then
convenient to represent $u$ in terms of a phase factor $z(t) = \exp (-i \varphi (t))$ and of an
amplitude $w(t)$ in such a way that asymptotically $\varphi (t)$ behaves as $\varphi^{(0)}(t)$ and
$w(t)$ tends to $w_+$. This is done by writing $u$ in the form \cite{7r} \cite{8r}
$$u(t) = M(t) \ D(t) \ z(t) \ w(t) \equiv \left ( \Lambda (w, \varphi ) \right ) (t) \quad .
\eqno(2.7)$$ 

\noi In I, we introduced three possible modified free evolutions $v_i(t)$ $i = 1,2,3$ and
correspondingly three parametrizations of $u(t)$ by $(w_i(t), \varphi_i(t))$, $i = 1,2,3$. The
choice (2.6) (2.7) corresponds to $i = 2$. We shall work exclusively with that choice throughout
this paper, and the subscript 2 is therefore consistently omitted. In I we used mostly the choice
$i = 3$ and dropped the subscript 3,  so that $(w, \varphi )$ in I means $(w_3, \varphi_3)$ as
opposed to $(w_2, \varphi_2)$ in this paper. This should be kept in mind when comparing results
from I and from this paper. \par

The construction of the wave operators for $u$ proceeds by first constructing the wave
o\-pe\-ra\-tors for the pair $(w , \varphi )$ and then recovering the wave operators for
$u$ therefrom by the use of (2.7). The evolution equation for $(w, \varphi )$ is obtained by substituting
(2.7) into the equation (1.1). One obtains the equation
$$\left ( i \partial_t + (2t^2)^{-1} \Delta - D^*\widetilde{g}D \right ) z w = 0 \eqno(2.8)$$

\noi for $zw$, with
$$\widetilde{g} \equiv \widetilde{g} \left ( |u|^2 \right ) = \widetilde{g} \left ( |Dw|^2 \right )
\quad , \eqno(2.9)$$

\noi or equivalently, by expanding the derivatives in (2.8),
$$\left \{ i \partial_t + (2t^2)^{-1} \Delta - i (2t^2)^{-1} \left ( 2 \nabla \varphi \cdot
\nabla + (\Delta \varphi ) \right ) \right \} w$$
$$+ \left \{ \partial_t \varphi - (2t^2)^{-1} \ |\nabla \varphi|^2 - D^*\widetilde{g}D \right \} w
= 0 \quad . \eqno(2.10)$$

We are now in the situation of a gauge theory. The equation (2.8) or (2.10) is invariant under the
gauge transformation $(w, \varphi ) \to (w \exp (i \sigma ), \varphi + \sigma )$, where $\sigma$ is
an arbitrary function of space time, and the original gauge invariant equation is not sufficient
to provide evolution equations for the two gauge dependent quantities $w$ and $\varphi$. At this
point we arbitrarily add the Hamilton-Jacobi equation as a gauge condition. This yields a system
of evolution equations for $(w, \varphi )$, namely

$$\hskip 2.5 truecm \left \{ \begin{array}{ll} \partial_t w = i(2t^2)^{-1} \Delta w + (2t^2)^{-1}
\left ( 2 \nabla \varphi \cdot \nabla + (\Delta \varphi ) \right ) w &\hskip 4 truecm (2.11) \\
& \\
\partial_t \varphi = (2t^2)^{-1} \ |\nabla \varphi |^2 + t^{- \gamma} \ g_0 (w, w)
&\hskip 4 truecm (2.12) \end{array} \right .$$

\noi where we have defined
$$ g_0(w_1, w_2) = \lambda \ {\rm Re} \ \omega^{\mu - n} \ w_1 \ \bar{w}_2 \eqno(2.13)$$

\noi and rewritten the nonlinear interaction term in (2.10) as
$$D^*\widetilde{g} \left ( |Dw|^2 \right ) D = t^{-\gamma} \ g_0(w, w) \quad . $$

\noi The gauge freedom in (2.11) (2.12) is now reduced to that given by an arbitrary function
of space only. It can be shown, actually it has been shown in I, that the Cauchy problem for
the system (2.11) (2.12) is locally wellposed in a neighborhood of infinity in time. The
solutions thereby obtained behave asymptotically as $w(t) = O(1)$ and $\varphi (t) = O(t^{1-
\gamma})$ as $t \to \infty$, a behaviour that is immediately seen to be compatible with
(2.11) (2.12). \par

We next study the asymptotic behaviour of the solutions of the auxiliary system (2.11) (2.12)
in more detail and try to construct wave operators for that system. For that purpose, we need
to choose a set of model functions playing the role of $v$, in the spirit of (2.6). In the
simple case $\gamma > 1/2$ considered in I, that set of model functions was taken to consist
of solutions of the system
$$\left \{ \begin{array}{l} \partial_t \ w^{(0)} = 0 \\ \\ \partial_t \ \varphi^{(0)} =
t^{-\gamma} \ g_0\left ( w^{(0)} , w^{(0)} \right ) \quad . \end{array} \right . \eqno(2.14)$$

\noi The general solution of (2.14) is
$$\left \{ \begin{array}{l} w^{(0)}(t) = w_+ \\ \\ \varphi^{(0)}(t) =
\psi_+ + \displaystyle{\int_1^t} dt_1 \ t_1^{- \gamma} \ g_0(w_+, w_+) \equiv \psi_+ + \varphi_0 (t) 
\end{array} \right . \eqno(2.15)$$

\noi and leads to (2.6) with $\varphi^{(0)} = \psi_+  + \varphi_0$. The asymptotic states for $(w,
\varphi )$ then consist of pairs $(w_+ , \psi_+)$. The choice (2.14) (2.15) is adequate for
$\gamma > 1/2$ because comparison of (2.11) (2.12) with (2.15) yields $\partial_t (\varphi -
\varphi_0) = O(t^{-2 \gamma})$ which is integrable at infinity for $\gamma > 1/2$, thereby allowing
for imposing an initial condition at $t = \infty$ for $\psi_1 = \varphi - \varphi_0$. For $\gamma
\leq 1/2$ however, the choice (2.14) (2.15) is not sufficient and one needs to construct more
accurate asymptotic functions. There are several ways to do that. The one we choose can be
motivated heuristically as follows. Let $p \geq 0$ be an integer. We write
$$\hskip 4 truecm \left \{ \begin{array}{ll} w = \displaystyle{\sum\limits_{0 \leq m \leq p}} w_m +
q_{p+1} \equiv W_p + q_{p+1} &\hskip 5 truecm (2.16) \\ & \\
\varphi = \displaystyle{\sum\limits_{0 \leq m \leq p}} \varphi_m + \psi_{p+1} \equiv \phi_p +
\psi_{p+1} &\hskip 5 truecm (2.17) \end{array} \right .$$

\noi with the understanding that asymptotically in $t$
$$w_m(t) = O \left ( t^{-m\gamma} \right ) \quad , \quad q_{p+1}(t) = o \left ( t^{-p\gamma} \right
) \quad , \eqno(2.18)$$
$$\varphi_m (t) = O \left ( t^{1 - (m+1)\gamma} \right ) \quad , \quad \psi_{p+1} (t) = o \left (
t^{1 - (p+1)\gamma} \right ) \quad . \eqno(2.19)$$

\noi Substituting (2.16) (2.17) into  (2.11) (2.12) and identifying the various powers of
$t^{-\gamma}$ yields the following system of equations for $(w_m, \varphi_m)$~: 
$$\hskip 1 truecm\left \{ \begin{array}{ll}\partial_t \ w_{m+1} = \left ( 2t^2 \right )^{-1}
\displaystyle{\sum\limits_{0 \leq j \leq m}} \left ( 2 \nabla \varphi_j \cdot \nabla + (\Delta
\varphi_j ) \right ) w_{m-j} &\hskip 2 truecm(2.20) \\  & \\ \partial_t \ \varphi_{m+1} = \left (
2t^2 \right )^{-1} \displaystyle{\sum\limits_{0 \leq j \leq m}} \nabla \varphi_j \cdot \nabla
\varphi_{m-j} + t^{- \gamma} \displaystyle{\sum\limits_{0 \leq j \leq m+1}} g_0 \left ( w_j ,
w_{m+1-j} \right ) &\hskip 2 truecm (2.21) \end{array}\right . $$

\noi for $m + 1 \geq 0$. Here it is understood that $w_j = 0$ and $\varphi_j = 0$ for $j < 0$, so
that the case $m = - 1$ of (2.20) (2.21) reduces to (2.14) with $w^{(0)} = w_0$ and
$\varphi^{(0)}  = \varphi_0$. We supplement that system with the initial conditions
$$\hskip 2.5 truecm \left \{ \begin{array}{ll} w_0(\infty ) = w_+ \qquad , \quad w_m(\infty ) = 0
\quad \hbox{for} \ m \geq 1 &\hskip 4.5 truecm (2.22) \\ & \\ \varphi_m (1) = 0 \qquad \hbox{for} \ 
0 \leq m \leq p \quad . &\hskip 4.5 truecm (2.23) \end{array} \right . $$

\noi The system (2.20) (2.21) with the initial conditions (2.22) (2.23) can be solved by
successive integrations~: knowing $(w_j, \varphi_j)$ for $0 \leq j \leq m$, one constructs
successively $w_{m+1}$ by integrating (2.20) between $t$ and $\infty$, and then $\varphi_{m+1}$
by integrating (2.21) between 1 and $t$. \par

If $(p+1)\gamma < 1$, that method of resolution reproduces the asymptotic behaviour in time (2.18)
(2.19) which was used in the first place to provide a heuristic derivation of the system (2.20)
(2.21). One can however consider that system and solve it by the same method for any integer $p$.
If $(p+1)\gamma > 1$, the asymptotic behaviour saturates at $w_m = O(t^{-1})$ for $m \gamma > 1$
and $\varphi_m = O(1)$ for $(m+1)\gamma > 1$. If $\gamma^{-1}$ is an integer, $(m + 1)\gamma = 1$
for some $m$, then $\varphi_m(t) =$\break \noindent $O$(Log $t$) and $w_{m+1} = O(t^{-1}$Log $t$).
\par

We now argue that for sufficiently large $p$, $\phi_p$ is a sufficiently good approximation for
$\varphi$ to ensure that $\psi_{p+1}$ has a limit as $t \to \infty$. In fact by comparing the
system (2.20) (2.21) with (2.11) (2.12), one finds that $\partial_t \ \psi_{p+1}$ is of the same
order in $t$ as $\partial_t \ \varphi_{p+1}$, namely $\partial_t \ \psi_{p+1} =
O(t^{-(p+2)\gamma})$, which is integrable at infinity for $(p + 2)\gamma > 1$. In this way every
solution $(w, \varphi )$ of the system (2.11) (2.12) as obtained previously has asymptotic states
consisting of $w_+ = \lim\limits_{t \to \infty} w(t)$ and $\psi_+ = \lim\limits_{t \to \infty}
\psi_{p+1}(t)$. \par

Conversely, under the condition $(p+2)\gamma > 1$, we shall be able to solve the system (2.11)
(2.12) by looking for solutions in the form (2.16) (2.17) with the additional initial condition
$\psi_{p+1} (\infty ) = \psi_+$, thereby getting a solution which is asymptotic to $(W_p, \phi_p
+ \psi_+)$ with 
$$w - W_p = O\left (t^{-(p+1)\gamma} \right ) \quad , \quad \varphi - \phi_p - \psi_+ = O    
\left ( t^{1-(p+2)\gamma} \right ) \quad . \eqno(2.24)$$

This allows to define a map $\Omega_0 : (w_+, \psi_+) \to (w, \varphi )$ which is essentially
the wave operator for $(w, \varphi )$. \par

It is an unfortunate feature of the methods used in this paper that both the construction of the
asymptotic states $(w_+, \psi_+ )$ of a given solution $(w, \varphi )$ and the construction of $(w,
\varphi )$ from given asymptotic states $(w_+, \psi_+)$ suffer from a loss of regularity of
roughly $p + 1$ derivatives, which prevents the two constructions to be inverse of each other in a
strict sense. \par

We next discuss the gauge covariance properties of $\Omega_0$. Two solutions $(w, \varphi )$ and
$(w' , \varphi ')$ of the system (2.11) (2.12) will be said to be gauge equivalent if they give
rise to the same $u$ through (2.7), namely if $w \exp (- i \varphi ) = w' \exp (- i \varphi ')$.
If $(w, \varphi )$ and $(w' , \varphi ')$ are two gauge equivalent solutions, one can show easily
that the difference $\varphi_- = \varphi ' - \varphi$ has a limit $\sigma$ when $t \to \infty$ and
that $w'_+ = w_+ \exp (i \sigma )$. Under that condition, it turns out that the phases
$\{\varphi_j\}$ and $\phi_p$ (but not the amplitudes) obtained by solving (2.20) (2.21) are
gauge invariant, namely $\varphi_m = \varphi '_m$ for $0 \leq m \leq p$ and therefore $\phi_p =
\phi '_p$, so that $\psi '_+ = \psi_+ + \sigma$. It is then natural to define gauge equivalence
of asymptotic states $(w_+, \psi_+)$ and $(w'_+ , \psi '_+)$ by the condition $w_+ \exp (- i
\psi_+) = w'_+ \exp (-i \psi '_+)$ and the previous result can be rephrased as the statement that
gauge equivalent solutions of (2.11) (2.12) in ${\cal R} (\Omega_0)$ have gauge equivalent
asymptotic states. Conversely, we are interested in showing that gauge equivalent asymptotic
states have gauge equivalent images under $\Omega_0$. Here however we meet with a technical
problem coming from the construction of $\Omega_0$ itself. For given $(w_+ , \psi_+)$ we
construct $(w, \varphi )$ in practice as follows. We take a (large) finite time $t_0$ and we
define a solution $(w_{t_0}, \varphi_{t_0})$ of the system (2.11) (2.12) by imposing a suitable
initial condition at $t_0$, depending on $(w_+, \psi_+)$, and using the known results for the
Cauchy problem with finite initial time. We then let $t_0$ tend to infinity and obtain $(w,
\varphi )$ as the limit of $(w_{t_0}, \varphi_{t_0})$. The simplest way to prove the gauge
equivalence of two solutions $(w, \varphi )$ and $(w', \varphi ')$ obtained in this way from
gauge equivalent $(w_+, \psi_+ )$ and $(w'_+, \psi '_+)$ consists in using an initial
condition at $t_0$ which already ensures that $(w_{t_0}, \varphi_{t_0})$ and $(w'_{t_0},
\varphi '_{t_0})$ are gauge equivalent. Unfortunately the natural choice $(w_{t_0}(t_0),
\varphi_{t_0} (t_0)) = (W_p(t_0), \phi_p (t_0) + \psi_+)$ does not satisfy that requirement as
soon as $p \geq 1$ because $\phi_p(t_0)$ is gauge invariant while $W_p (t_0) \exp (- \psi_+ )$
is not. In order to overcome that difficulty, we introduce a new amplitude $V$ and a new phase
$\chi$ defined by solving the transport equations
$$\hskip 4 truecm \left \{ \begin{array}{ll} \partial_t V = (2t^2)^{-1} \left ( 2 \nabla \phi_{p-1}
\cdot \nabla + \left ( \Delta \phi_{p-1} \right ) \right ) V &\hskip 3.8 truecm(2.25) \\ & \\
\partial_t \chi = t^{-2} \ \nabla \phi_{p-1} \cdot \nabla \chi &\hskip 3.8 truecm(2.26) \end{array}
\right . $$

\noi with initial condition
$$V(\infty ) = w_+ \qquad , \qquad \chi (\infty ) = \psi_+ \quad . \eqno(2.27)$$

\noi It follows from (2.25) (2.26) that $V \exp (- i \chi )$ satisfies the same transport
equation as $V$, now with gauge invariant initial condition $(V \exp (- i \chi ))(\infty ) =
w_+ \exp (-i \psi_+)$, and is therefore gauge invariant. Furthermore, $(V, \chi )$ is a
sufficiently good approximation of $(W_p , \psi_+)$ in the sense that
$$V(t) - W_p(t) = O \left (t^{-(p+1)\gamma} \right ) \quad , \quad \chi (t) - \psi_+ =
O(t^{-\gamma}) \quad . \eqno(2.28)$$

\noi One then takes $(w_{t_0} (t_0), \varphi_{t_0} (t_0)) = (V(t_0), \phi_p (t_0) + \chi
(t_0))$ as an initial condition at time $t_0$, thereby ensuring that $(w_{t_0},
\varphi_{t_0})$ and $(w'_{t_0}, \varphi '_{t_0})$ are gauge equivalent. That equivalence is
easily seen to be preserved in the limit $t_0 \to \infty$. Furthermore, the estimates (2.28)
ensure that the asymptotic properties (2.24) are preserved by the modified construction. As a
consequence of the previous discussion, the map $\Omega_0$ is gauge covariant, namely induces
an injective map of gauge equivalence classes of asymptotic states $(w_+, \psi_+)$ to gauge
equivalence classes of solutions $(w, \varphi )$ of the system (2.11) (2.12). \par

The wave operator for $u$ is obtained from $\Omega_0$ just defined and from $\Lambda$ defined by
(2.7). From the previous discussion it follows that the map $\Lambda \circ \Omega_0 : (w_+, \psi_+ )
\to u$ is injective from gauge equivalence classes of asymptotic states $(w_+ , \psi_+)$ to
solutions of (1.1). In order to define a wave operator for $u$ involving only the asymptotic state
$u_+$ but not an arbitrary phase $\psi_+$, we choose a representative in each equivalence class
$(w_+, \psi_+)$, namely we define the wave operator for $u$ as the map $\Omega : u_+ \to u =
(\Lambda \circ \Omega_0) (Fu_+, 0)$. Since each equivalence class of asymptotic states contains at
most one element with $\psi_+ = 0$, the map $\Omega$ is again injective. We shall prove in
addition that ${\cal R}(\Omega ) = {\cal R} (\Lambda \circ \Omega_0)$ if $p \leq 2$. (This need not
be the case if $p \geq 3$, because derivative losses in the construction generate a mismatch
between the regularity properties required on $w_+$ and $\psi_+$, so that gauge equivalence
classes of asymptotic states need not contain an element with $\psi_+ = 0$ in that case). \par

The previous heuristic discussion was based in part on a number of asymptotic estimates in terms
of negative powers of $t$. However if $\gamma^{-1}$ is an integer some of these estimates have
to be replaced or supplemented by logarithms. In order to treat all values of $\gamma \in
(0,1]$ in a unified way, we shall introduce a number of estimating functions of time defined by
integral representations. Those functions are smooth in $\gamma$, in particular at integer
values of $\gamma^{-1}$. They generate the logarithms automatically whenever needed, and they
recombine nicely between themselves in the derivation of the main estimates. The simplest
example thereof is $h_0(t)$ defined by (3.19) below. \par

In the same way as in I, the system (2.11) (2.12) can be rewritten as a system of equations for $w$
and for $s = \nabla \varphi$, from which $\varphi$ can then  be recovered by (2.12), thereby
leading to a slightly more general theory since the system for $(w, s)$ can be studied without
even assuming that $s$ is a gradient. In I, we first studied the system for $(w, s)$ and then
deduced therefrom the relevant results for $(w , \varphi )$. Here for simplicity we shall use
exclusively the variables $(w, \varphi )$. The same remark applies to the system (2.20) (2.21).
\par

We are now in a position to describe in more detail the contents of the technical parts of this
paper, namely Sections 3-7. In Section 3, we introduce some notation, we define the relevant
function spaces needed to study the system (2.11) (2.12), we recall from I a number of Sobolev
and energy estimates, we then introduce the estimating functions of time mentioned above and we
derive a number of estimates for them. In Section 4, we recall from I some preliminary results
on the Cauchy problem for the auxiliary system (2.11) (2.12) and on the asymptotic behaviour of
its solutions. In Section 5 we study the systems (2.20) (2.21) and (2.25) (2.26) defining the
asymptotic dynamics. We first derive a number of properties and estimates for the solutions of
the system (2.20) (2.21), defined inductively (Proposition 5.1). We then prove the existence
and some properties of the solutions of the transport equations (2.25) and (2.26) 
in a slightly more general setting (Propositions 5.2 and 5.3 respectively). We finally
specialize those results to the case at hand and compare $V$ with $W_p$ defined by (2.16)
(Proposition 5.4). In Section 6 we study in detail the asymptotic behaviour in time of
solutions of the auxiliary system (2.11) (2.12). We first derive asymptotic estimates on the
approximation of the available solutions $(w , \varphi )$ of that system by the asymptotic
functions $(W_m , \phi_m)$ defined by (2.16) (2.17), and in particular we complete the proof
of existence of asymptotic states for those solutions (Proposition 6.1). We then turn to the
construction of local wave operators at infinity. For a given solution $(V, \chi )$ of the
system (2.25) (2.26) and a given (large) $t_0$, we construct a solution $(w_{t_0},
\varphi_{t_0})$ of the system (2.11) (2.12) which coincides with $(V, \phi_p + \chi )$ at
$t_0$ and we estimate it uniformly in $t_0$ (Proposition 6.2). We then prove that when $t_0
\to \infty$, $(w_{t_0}, \varphi_{t_0})$ has a limit $(w, \varphi )$ which is asymptotic both
to $(V, \phi_p + \chi )$ and to $(W_p , \phi_p + \psi_+)$ (Proposition 6.3). Finally in
Section 7, we exploit the results of Section 6 to construct the wave operators for the
equation (1.1) and to describe the asymptotic behaviour of solutions in their range. We
first prove that the local wave operator at infinity for the system (2.11) (2.12) defined
through Proposition 6.3 in Definition 7.1 is gauge covariant in the sense of Definitions 7.2
and 7.3 in the best form that can be expected with the available regularity (Propositions
7.2 and 7.3). With the help of some information on the Cauchy problem for (1.1) at finite
time (Proposition 7.1), we then define the wave operator $\Omega : u_+ \to u$ (Definition
7.4), we prove that it is injective and under suitable restrictions, that it has the
expected range (Proposition 7.4). We then collect all the available information on $\Omega$
and on solutions of (1.1) in its range in Proposition 7.5, which contains the main results
of this paper.

\section{Notation and preliminary estimates}
\hspace*{\parindent} In this section, we define the function spaces where we shall study the
auxiliary system (2.11) (2.12) and we recall from I a number of Sobolev and energy type estimates
which hold in those spaces. We then introduce a number of estimating functions of time and we derive
a number of relations and estimates for them. \par

We shall use Sobolev spaces of integer order $H_r^k$ defined for $1 \leq r \leq \infty$ by
$$H_r^k = \left \{ u : \parallel u; H_r^k \parallel \ \equiv \sum_{0 \leq j \leq k} \parallel
\partial^j u \parallel_r \ < \infty \right \}$$

\noi and the associated homogeneous spaces $\dot{H}_r^k$ with norm
$$\parallel u; \dot{H}_r^k \parallel \ = \ \parallel \partial^k u \parallel_r$$

\noi where
$$\parallel \partial^j u \parallel_r \ = \sum_{\alpha : |\alpha | = j} \parallel \partial^{\alpha}
u \parallel_r \quad .$$

\noi The subscript $r$ will be omitted if $r = 2$. \par

Let $\ell_0 = [n/2]$ and define $r_0$ by $\delta (r_0) = \ell_0$ so that $r_0 = 2n$ for odd $n$
and $r_0 = \infty$ for even $n$. Let $k$ and $\ell$ be nonnegative integers with $\ell \geq \ell_0
- 1$. We shall look for $w$ as a complex valued function in spaces $L_{loc}^{\infty}(I, H^k)$ or
${\cal C}(I, H^k)$ and for $\varphi$ as a real valued function in spaces $L_{loc}^{\infty}(I,
Y^{\ell})$ or ${\cal C}(I, Y^{\ell})$ where 
$$Y^{\ell} = L^{\infty} \cap \dot{H}_{r_0}^1 \cap \dot{H}^{\ell_0 + 1} \cap \dot{H}^{\ell + 2}
\quad . \eqno(3.1)$$

\noi The spaces $Y^{\ell}$ are easily seen to be duals of Banach spaces and satisfy the embedding
$Y^{\ell '} \subset Y^{\ell}$ for $\ell ' \geq \ell$. We shall use systematically the notation
$$|w|_k = \ \parallel w; H^k\parallel \qquad , \qquad |\varphi |_{\ell} = \ \parallel \varphi ;
Y^{\ell} \parallel \eqno(3.2)$$

\noi and the meaning of the symbol $|a|_b$ will be made unambiguous by the fact that the pair $(a,
b)$ contains either the pair $(w, k)$ or the pair $(\varphi , \ell )$. Note that the second
notation in (3.2) is different from, although closely related to, the similar notation in I
which was used for $s = \nabla \varphi$. \par

We recall the following result from I (see Lemma I.3.5). \\

\noi {\bf Lemma 3.1.} {\it Let $\varphi$ be a real function with $\nabla \varphi \in L^{\infty}
\cap \dot{H}^{\ell}$ for some $\ell > n/2$ and let $k \leq \ell + 1$. Then the following estimate
holds~:}
$$\left | \exp (- i \varphi ) w \right |_k \leq C \left ( 1 + \ \parallel \nabla \varphi ;
L^{\infty} \cap \dot{H}^{\ell} \parallel \right )^k \ |w|_k \quad . \eqno(3.3)$$

\noi {\it Let in addition $\varphi \in L^{\infty}$. Then the following estimate holds~:}
$$\left | \left ( \exp (- i \varphi ) - 1 \right ) w \right |_k \leq C \left (\parallel
\varphi \parallel_{\infty}  \ + \ \parallel \nabla \varphi ; L^{\infty} \cap \dot{H}^{\ell} \parallel
\left ( 1 + \ \parallel \nabla \varphi ; L^{\infty} \cap \dot{H}^{\ell} \parallel \right )^{k-1}
\right ) |w|_k  \ .\eqno(3.4)$$ \vskip 3 truemm

In order to state the estimates that are relevant for the study of the system (2.11) (2.12), it is
useful to give the following definition (see Definition I.3.1). \\

\noi {\bf Definition 3.1.} Let $0 < \mu \leq n - 2$. A pair of nonnegative integers $(k, \ell )$
will be called admissible if it satisfies $k \leq \ell$, $\ell > n/2$ and
$$\ell + 2 + \mu \leq (n/2 + 2k) \wedge (n + k)  \eqno(3.5)$$

\noi and in addition $k > n/2$ if $\ell + 2 + \mu = n + k$ and
$$n/2 + 3 + \mu < (n/2 + 2k) \wedge (n + k)$$

\noi if $n$ is even. \\

For $\mu = n - 2$, admissible pairs are pairs $(k, \ell )$ such that $k = \ell > n/2$. If $(k, \ell
)$ is admissible, so is $(k + j, \ell + j)$ for any positive integer $j$. Admissible pairs always
have $k \geq 2$. For $n = 3$, $\mu = 1$, the pair (2,2) is admissible. \par

The following Sobolev like inequalities will be essential to study the system (2.11) (2.12). \\

\noi {\bf Lemma 3.2.} {\it Let $\ell > n/2$ and $k \leq \ell$. Then the following estimates hold~:}
$$\left | \left ( 2 \nabla \varphi \cdot \nabla + (\Delta \varphi ) \right ) w\right |_{k-1} \leq C
 | \varphi |_{\ell - 1} \ |w|_k \quad , \eqno(3.6)$$
$$\left | \nabla \varphi_1 \cdot \nabla \varphi_2 \right |_{\ell - 1} \leq C |\varphi_1|_{\ell} \
|\varphi_2|_{\ell} \quad . \eqno(3.7)$$

\noi {\it Assume in addition that $(k , \ell )$ is admissible. Let $g_0$ be defined by (2.13). Then}
$$\left | g_0 (w_1 \ w_2) \right |_{\ell } \leq C |w_1|_k \ |w_2|_k \quad , \eqno(3.8)$$
$$\left | g_0 (w_1 \ w_2) \right |_{\ell - 1} \leq C |w_1|_k \ |w_2|_{k-1} \quad . \eqno(3.9)$$
\vskip 3 truemm

\noi {\bf Sketch of proof.} (3.6) follows from Lemma I.3.4 by the same estimates as in Lemma I.3.9.
The estimate (3.7) essentially follows from Lemma I.3.3. The estimates (3.8) and (3.9) follow from
Corollary I.3.1. \par \nobreak
\hfill $\sq$ \par

In addition to the previous estimates, we shall need energy type estimates for solutions of the
following transport equations
$$\partial_t w  = (2t^2)^{-1} \Big \{ i \theta \Delta w + \left ( 2 \nabla \phi \cdot \nabla +
(\Delta \phi ) \right ) w \Big \} + R_1 \quad , \eqno(3.10)$$
$$\partial_t \varphi = (2t^2)^{-1} \Big \{ \theta |\nabla \varphi |^2 + 2 \nabla \phi \cdot
\nabla \varphi \Big \} + R_2 \quad  \eqno(3.11)$$

\noi where $\theta$ is a real constant and $\phi$, $R_1$, $R_2$ are given functions of space time.
Those estimates will be stated in differential form for brevity, although they should be understood
in integrated from. They hold for functions that are sufficiently regular in time, for instance
locally bounded in the relevant norms. \\

\noi {\bf Lemma 3.3.} {\it Let $\ell > n/2$ and $k \leq \ell$. \par
(1) Let $w$ satisfy (3.10). Then the following estimate holds~:}
$$\left | \partial_t|w|_k \right | \leq C \ t^{-2} |\phi|_{\ell} \ |w|_k + |R_1|_k \quad .
\eqno(3.12)$$

{\it (2) Let $\varphi$ satisfy (3.11). Then the following estimates hold~:}
$$\left | \partial_t |\varphi|_{\ell} \right | \leq C \ t^{-2} | \varphi |_{\ell} \left ( |\theta |
\ |\varphi|_{\ell}  + |\phi |_{\ell + 1} \right ) + |R_2|_{\ell} \quad , \eqno(3.13)$$  
$$\left | \partial_t |\varphi|_{\ell - 1} \right | \leq C \ t^{-2} | \varphi |_{\ell - 1} \left (
|\theta | \ |\varphi|_{\ell}  + |\phi |_{\ell} \right ) + |R_2|_{\ell - 1} \quad . \eqno(3.14)$$

\noi {\bf Sketch of proof.} \par
(3.12) follows from Lemmas I.3.2 and I.3.4 by the same estimates as in Lemma I.3.7. \par

(3.13) and (3.14) follow from Lemmas I.3.2 and I.3.3 by the same estimates as in Lemmas I.3.7 and
I.3.9. \par \nobreak
\hfill $\sq$ \par

\noi {\bf Lemma 3.4.} {\it Let $\ell > n/2$ and $k \leq \ell$. Let $w$ and $\varphi$ satisfy
(3.10) and (3.11) respectively, with $\theta = 0$, $R_1 = 0$ and $R_2 = 0$. Then the following
estimates hold~:}
$$\left | \partial_t |w|_{k+1} \right | \leq C \ t^{-2} \left ( |\phi |_{\ell} \ |w|_{k+1} + |\phi
|_{\ell + 1} \ |w|_k \right ) \quad ,  \eqno(3.15)$$  
$$\left | \partial_t |\varphi |_{\ell +1} \right | \leq C \ t^{-2} \left ( |\phi |_{\ell} \
|\varphi|_{\ell +1} + |\phi |_{\ell + 2} \ |\varphi|_{\ell} \right ) \quad .  \eqno(3.16)$$
\vskip 3 truemm

\noi {\bf Sketch of proof.} \par
(3.15) and (3.16) follow from Lemmas I.3.2, I.3.3 and I.3.4 by the same estimates as in Lemma
I.3.8. \par \nobreak
\hfill $\sq$   \par

\noi {\bf Lemma 3.5.} {\it Let $\ell > n/2$ and $k \leq \ell$. Let $w_1$, $w_2$ and $\varphi_1$,
$\varphi_2$ satisfy (3.10) and (3.11) with $\phi = \phi_1$ and $\phi = \phi_2$ respectively, and
with $\theta = 0$, $R_1 = 0$ and $R_2 = 0$. Let $w_- = w_1 - w_2$, $\varphi_- = \varphi_1 -
\varphi_2$ and $\phi_- = \phi_1 - \phi_2$. Then the following estimates hold~:}
$$\left | \partial_t |w_-|_k \right | \leq C \ t^{-2} \Big ( |\phi_2|_{\ell} \ |w_-|_k +
|\phi_-|_{\ell} \ |w_1|_k + \parallel \nabla \phi_- \parallel_{\infty} \ |w_1|_{k+1} \Big )
\quad , \eqno(3.17)$$  
$$\left | \partial_t |\varphi_-|_{\ell} \right | \leq C \ t^{-2} \Big ( |\phi_2|_{\ell +1} \
|\varphi_-|_{\ell} + |\phi_-|_{\ell + 1} \ |\varphi_1|_{\ell} \ + \parallel \nabla
\phi_- \parallel_{\infty} \ |\varphi_1|_{\ell +1} \Big ) \quad . \eqno(3.18)$$
\vskip 3 truemm

\noi {\bf Sketch of proof.} \par
(3.17) and (3.18) follow from Lemmas I.3.2, I.3.3 and I.3.4 by the same estimates as in Lemma
I.3.10. \par \nobreak
\hfill $\sq$ \par

We now introduce a number of estimating functions of time and derive a number of estimates and
relations for them. We start with 
$$h_0 (t) = \int_1^t dt_1 \ t_1^{- \gamma} \eqno(3.19)$$

\noi so that
$$h_0 (t) = \left \{ \begin{array}{ll} (1 - \gamma )^{-1} (t^{1- \gamma} - 1) &\qquad
\hbox{for} \ \gamma \not= 1 \\ & \\ {\rm Log} \ t &\qquad \hbox{for} \ \gamma = 1 \quad .\end{array}
\right . \eqno(3.20)$$

\noi The basic building block for the subsequent functions is the function $h$ defined by
$$h(t) = \int_1^{\infty} dt_1 \ t_1^{-\gamma } (t \vee t_1)^{-1} \quad , \eqno(3.21)$$
\noi which can also be written as
$$h(t) = t^{-1} \ h_0 (t) + \gamma^{-1} \ t^{-\gamma} = \int_t^{\infty} dt_1 \ t_1^{-2} \ h_0 (t_1)
\eqno(3.22)$$
\noi and is explicitly computed as

$$h (t) = \left \{ \begin{array}{ll} \gamma^{-1} (1 - \gamma )^{-1} (t^{- \gamma} - \gamma t^{-1})
&\qquad \hbox{for} \ \gamma \not= 1 \\ &\\ t^{-1} (1 + {\rm Log} \ t) &\qquad \hbox{for} \ \gamma = 1
\quad .\end{array} \right . \eqno(3.23)$$

\noi It follows from (3.21) that $t \ h(t)$ is increasing in $t$ and from (3.20) (3.23) that $t \
h(t) \ h_0(t)^{-1}$ is decreasing in $t$. The function $h$ satisfies the estimates
$$\gamma^{-1} \left ( t^{-\gamma} \vee t^{-1} \right ) \leq h(t) \leq |1 - \gamma|^{-1} \left (
\gamma^{-1}  \ t^{-\gamma} \vee t^{-1} \right ) \quad . \eqno(3.24)$$

\noi The first inequality in (3.24) follows in part from (3.22) and in part from the monotony
of\break \noindent $t$ $h(t)$, while the second inequality follows from (3.23) and holds only for
$\gamma \not= 1$. \par

In some cases where we shall need to indicate the dependence of $h_0$ and $h$ on $\gamma$, we shall
write $h_0 (\gamma , t)$ and $h(\gamma , t)$ for $h_0(t)$ and $h(t)$. \par

We next define for any $m \geq 0$
$$N_m (t) = \int_1^t dt_1 \ t_1^{-\gamma} \ h^m(t_1) \quad , \eqno(3.25)$$
$$Q_m(t) = \int_1^{\infty} dt_1 \ t_1^{- \gamma} (t \vee t_1)^{-1} \ h^m(t_1) \quad , \eqno(3.26)$$

\noi so that $N_0 = h_0$ and $Q_0 = h$. Those functions are smooth in $\gamma$. Clearly $N_m$ is
increasing and $Q_m$ is decreasing in $t$, while $t \ Q_m(t)$ is increasing in $t$, so that $Q_m(t)
\geq Q_m(1) \ t^{-1}$. From the fact that $h$ is decreasing, it follows that 
$$N_{i+j}(t) \leq \gamma^{-j} \ N_i(t) \leq \gamma^{-(i+j)} \ h_0(t) \quad , \eqno(3.27)$$
$$Q_{i+j}(t) \leq \gamma^{-j} \ Q_i(t) \leq \gamma^{-(i+j)} \ h(t) \eqno(3.28)$$

\noi for all $i \geq 0$, $j \geq 0$. It follows from (3.24) that $N_m$ and $Q_m$ satisfy the lower
and upper bounds
$$\gamma^{-m} \ h_0 \left ( (m + 1) \gamma , t \right ) \leq N_m (t) \leq (1 - \gamma)^{-m} \
\gamma^{-m} \ h_0 \left ( (m + 1) \gamma , t \right ) \eqno(3.29)$$
$$\gamma^{-m} \ h \left ( (m + 1) \gamma , t \right ) \leq Q_m (t) \leq (1 - \gamma)^{-m} \
\gamma^{-m} \ h \left ( (m + 1) \gamma , t \right ) \eqno(3.30)$$

\noi where the lower bounds hold for all $\gamma > 0$ and the upper bounds for $0 < \gamma < 1$ if
$m \geq 1$. From (3.20) and (3.23), it follows that $N_m(t)$ and $Q_m(t)$ behave as
$t^{1-(m+1)\gamma}$ and $t^{-(m+1)\gamma}$ respectively as $t \to \infty$ if $(m + 1) \gamma < 1$.
If $(m + 1) \gamma = 1$, $N_m(t)$ and $Q_m(t)$ produce logarithms and behave as Log $t$ and
$t^{-1}$ Log $t$ respectively as $t \to \infty$. If $(m + 1) \gamma > 1$, $N_m(t)$ and $Q_m(t)$
saturate respectively as Constant and $t^{-1}$ when $t \to \infty$. For $m \geq 1$, the upper
bounds in (3.29) and (3.30) blow up when $\gamma$ tends to one, but the same conclusions still
hold. \par

For $(m + 2) \gamma > 1$, we finally define
$$P_m (t) = \int_1^{\infty} dt_1 \ t_1^{-\gamma} \ h(t \vee t_1) \ h^m(t_1) \quad , \eqno(3.31)$$
$$R_m(t) = \int_t^{\infty} dt_1 \ t_1^{-2} \ P_m (t_1) \quad . \eqno(3.32)$$

In particular $P_0$ is explicitly computed as
$$P_0(t) = h_0(t) \left ( h(t) + t^{-1} \ t^{-\gamma} \right ) + 2 \gamma^{-1}(2 \gamma - 1)^{-1} \
t^{1 - 2 \gamma} \quad . \eqno(3.33)$$

\noi Clearly $P_m(t)$ and $R_m(t)$ are decreasing in $t$, while $P_m(t) \ h(t)^{-1}$ is increasing
in $t$, so that $P_m(t) \geq P_m(1) \gamma h(t)$. It follows from (3.24) that $P_m$ satisfies the
lower and upper bounds 
$$P_m(t) \left \{ \begin{array}{l} \geq 1 \\ \leq (1 - \gamma )^{-(m+1)} \end{array} \right \}
\gamma^{-(m+1)} \left ( t^{-\gamma} \ h_0 \left ( (m+1) \gamma , t \right ) + \left ( (m + 2) \gamma
- 1 \right )^{-1} \ t^{1-(m+2)\gamma} \right ) \quad .  \eqno(3.34)$$

\noi From (3.20) it follows that $P_m(t)$ behaves as $t^{1-(m+2)\gamma}$ as $t \to \infty$ if $(m +
1)\gamma < 1$. If $(m + 1) \gamma = 1$, $P_m (t)$ behaves as $t^{-\gamma}$ Log $t$. If $(m + 1)
\gamma > 1$, $P_m(t)$ saturates at $t^{-\gamma}$ as long as $\gamma < 1$. \par

We now collect a number of relations and estimates satisfied by the previous estimating functions.
\\

\noi {\bf Lemma 3.6.} {\it Let $i$, $j$ and $m$ be nonnegative integers. Let $1 \leq a \leq b$ and
$t \geq 1$. Then the following identities and estimates hold~:} 
$$\int_t^{\infty} dt_1 \ t_1^{-2} \ N_m(t_1) = Q_m(t) \eqno(3.35)$$
$$\int_1^t dt_1 \ t_1^{-2} \ h_0(t_1) \ N_m(t_1) = N_{m+1}(t) - h(t) \ N_m(t) \leq N_{m+1}(t)
\eqno(3.36)$$
$$\int_t^{\infty} dt_1 \ t_1^{-2} \ h_0(t_1) \ N_m(t_1) = P_m(t) \qquad \hbox{for} \ (m+2)\gamma >
1 \eqno(3.37)$$
$$\int_a^b dt \ t^{-2} \ N_i(t) \ N_j(t) \leq \int_a^b dt \ t^{-2} \ h_0(t) \ N_{i+j}(t)
\eqno(3.38)$$
$$\int_a^b dt \ t^{-2} \ N_i(t) \ Q_j(t) \leq \int_a^b dt \ t^{-2} \ h(t) \ N_{i+j}(t) \leq
\int_a^b dt \ t^{-2} \ N_{i+j+1} (t) \eqno(3.39)$$
$$\int_a^b dt \ t^{-\gamma} \ Q_i (t) \ Q_j(t) \leq \int_a^b dt \ t^{-\gamma} \ h(t) \ Q_{i+j}(t)
\eqno(3.40)$$
$$\int_t^{\infty} dt_1 \ t_1^{-\gamma} \ h(t_1) \ Q_{m-1}(t_1) \leq \int_t^{\infty} dt_1 \
t_1^{-\gamma} \ Q_m(t_1) \qquad \hbox{for} \ m\geq 1 \ , \ (m + 2)\gamma > 1 \quad . \eqno(3.41)$$
$$\int_t^{\infty} dt_1 \ t_1^{-\gamma} \ Q_m (t_1) \leq P_m(t) \qquad \hbox{for} \ (m+2) \gamma >
1 \quad . \eqno(3.42)$$
$$\int_1^t dt_1 \ t_1^{-\gamma} \ h(t_1) \ Q_{m-1}(t_1) \leq N_{m+1}(t) \eqno(3.43)$$
$$\int_1^t dt_1 \ t_1^{-\gamma} \ Q_m(t_1) \leq N_{m+1}(t) \eqno(3.44)$$
$$\int_a^b dt \ t^{-\gamma} \ Q_m(t) \leq Q_m(a) \left ( h_0 (b) - h_0(a) \right ) \eqno(3.45)$$
$$\int_a^b dt \ t^{-\gamma} \ h(t) \ Q_{m-1}(t) \leq 2 Q_m (a) \left ( h_0(b) - h_0(a) \right )
\eqno(3.46)$$
$$R_m(t) \leq C_m \ h(t) \ Q_m(t) \qquad \hbox{for} \ (m + 2)\gamma > 1 \eqno(3.47)$$

\noi where
$$C_m = (2m + 3) \gamma \left ( (m + 2) \gamma - 1 \right )^{-1} \quad .$$ 
\vskip 3 truemm

\noi {\bf Proof.} \par
(3.35). By the definition of $N_m$ and $Q_m$
$$\int_t^{\infty} dt_1 \ t_1^{-2} \ N_m(t_1) = \int_t^{\infty}dt_1 \ t_1^{-2} \int_1^{t_1} dt_2 \
t_2^{-\gamma} \ h^m(t_2)$$
$$= \int_1^{\infty} dt_2 \ t_2^{-\gamma} \ h^m(t_2) \int_{t \vee t_2}^{\infty} dt_1 \ t_1^{-2} =
\int_1^{\infty} dt_2 \ t_2^{-\gamma} \left ( t \vee t_2 \right )^{-1} \ h^m(t_2) = Q_m(t) \quad .$$

(3.36). By the definition of $N_m$ and integration by parts
$$\int_1^t dt_1 \ t_1^{-2} \ h_0(t_1) \ N_m(t_1) = - \int_1^t dt_1 \ h'(t_1) \ N_m(t_1)$$
$$= - h(t) \ N_m(t) + \int_1^t h(t) \ N'_m(t) = N_{m+1}(t) - h(t) \ N_m(t) \quad .$$

(3.37). By the definitions of $N_m$ and $P_m$ and integration by parts
$$\int_t^{\infty} dt_1 \ t_1^{-2} \ h_0(t_1) \ N_m(t_1) = h(t) \ N_m(t) + \int_t^{\infty} dt_1 \
t_1^{-\gamma} \ h^{m+1}(t_1)$$
$$= \int_1^{\infty} dt_1 \ t_1^{-\gamma} \ h\left ( t \vee t_1 \right ) \ h^m(t_1) = P_m(t) \quad .$$

(3.38). By the definition of $N_m$
$$\int_a^b dt \ t^{-2} \ N_i(t) \ N_j(t) = \int_a^b dt \ t^{-2} \int_1^t dt_1 \ t_1^{-\gamma} \ h^i
(t_1) \int_1^t dt_2 \ t_2^{-\gamma} \ h^j (t_2) \quad .$$

\noi For fixed $i + j$, the last integral is logarithmically convex in $i$ (or $j$) and therefore
estimated by the maximum of its values for $i = 0$ and $j = 0$, which are equal by symmetry and
equal to the RHS of (3.38). \par \vskip 3 truemm

(3.39). By the definition of $N_m$ and $Q_m$ 
$$\int_a^b dt \ t^{-2} \ N_i(t) \ Q_j(t) = \int_a^b dt \ t^{-2} \int_1^t dt_1 \ t_1^{-\gamma} \
h^i(t_1) \int_1^{\infty} dt_2 \ t_2^{-\gamma} \left ( t \vee t_2 \right )^{-1} \ h^j (t_2) \quad .$$

\noi We split the integral over $t_2$ into the subregions $t_2 \leq t$ and $t_2 \geq t$. In the
region $t_2 \leq t$, by logarithmic convexity and symmetry, we estimate the integral by replacing
$h ^i(t_1)$ $h^j(t_2)$ by $h^{i+j}(t_1)$. In the region $t_2 \geq t$, we make the same
replacement because $t_2 \geq t \geq t_1$ and $h$ is decreasing in $t$. We obtain 
$$\cdots \leq \int_a^b dt \ t^{-2} \int_1^t dt_1 \ t_1^{-\gamma} \ h^{i+j}(t_1) \
h(t) =  \int_a^b dt \ t^{-2} \ h(t) \ N_{i+j}(t)$$

\noi which yields the first inequality in (3.39). Using in addition the fact that $h(t) \leq
h(t_1)$ for $t_1 \leq t$ yields the second inequality. \par \vskip 3 truemm

(3.40). By the definition of $Q_m$, the LHS of (3.40) is logarithmically convex in $i$ or $j$ for
fixed $i + j$, and symmetric in $i$ and $j$, and is therefore estimated by its end point values,
namely with $i$, $j$ replaced by 0 and $i + j$. \par \vskip 3 truemm

(3.41) and (3.42). By the definition of $Q_m$
$$\int_t^{\infty} dt_1 \ t_1^{-\gamma} \ h(t_1) \ Q_{m-1} (t_1) = \int_t^{\infty} dt_1 \
t_1^{-\gamma} \ h(t_1) \int_1^{\infty} dt_2 \ t_2^{-\gamma} \left ( t_1 \vee t_2 \right )^{-1} \
h^{m-1} (t_2) \quad . \eqno(3.48)$$

\noi We estimate the last integral by replacing $h(t_1)$ by $h(t_2)$, by logarithmic convexity and
symmetry in the region $t_2 \geq t$ and by monotony of $h$ in the region $t_2 \leq t (\leq t_1)$,
thereby continuing (3.48) by 
$$\cdots \leq \int_t^{\infty} dt_1 \ t_1^{-\gamma} \int_1^{\infty} dt_2 \ t_2^{-\gamma} \left (
t_1 \vee t_2 \right )^{-1} \ h^m(t_2)$$

\noi which is the RHS of (3.41) and the LHS of (3.42), 
$$\cdots = \int_1^{\infty} dt_2 \ t_2^{-\gamma} \ h^m(t_2) \ \int_t^{\infty} dt_1 \ t_1^{-\gamma}
\left ( t_1 \vee t_2 \right )^{-1} \quad .$$

\noi We estimate the last integral by $h(t \vee t_2)$ by first replacing $t_1 \vee t_2$ by $t_1
\vee t \vee t_2$, since $t_1 \geq t$, and then extending the integration over $t_1$ to $[1, \infty
)$, thereby obtaining
$$\cdots \leq \int_1^{\infty} dt_2 \ t_2^{-\gamma} \ h^m(t_2) \ h(t \vee t_2) = P_m(t) \quad .$$

(3.43) and (3.44). By the definition of $Q_m$
$$\int_1^t dt_1  \ t_1^{-\gamma} \left ( h(t_1) Q_{m-1}(t_1) \ \hbox{or} \ Q_m(t_1) \right )$$
$$= \int_1^t dt_1 \ t_1^{-\gamma} \int_1^{\infty} dt_2 \ t_2^{-\gamma} \left ( t_1 \vee t_2 \right
)^{-1} \left ( h(t_1) \ h^{m-1}(t_2) \ \hbox{or} \ h^m(t_2) \right )$$
$$\leq \int_1^t dt_1 \ t_1^{-\gamma} \ h^m(t_1) \int_1^{\infty} dt_2 \ t_2^{-\gamma} \left ( t_1
\vee t_2 \right )^{-1} = N_{m+1}(t)$$

\noi by logarithmic convexity and symmetry in the region $t_2 \leq t$ and by monotony of $h$ in the
region $t_2 \geq t(\geq t_1)$. \par \vskip 3 truemm

(3.45) follows immediately from the fact that $Q_m$ is decreasing in $t$. \par \vskip 3 truemm

(3.46). We first prove that 
$$h(t) \ Q_{m-1}(t) \leq 2Q_m(t) \quad . \eqno(3.49)$$

\noi In fact
$$Q_{m-1}(t) \  h(t) = \int_1^{\infty} dt_1 \ t_1^{-\gamma} \left ( t \vee t_1 \right )^{-1}
\int_1^{\infty} dt_2 \ t_2^{-\gamma} \left ( t \vee t_2 \right )^{-1} \ h^{m-1} (t_1)$$
$$= \int_1^{\infty} dt_1 \ t_1^{-\gamma} \left ( t \vee t_1 \right )^{-1} \int_{t_1}^{\infty} dt_2
\ t_2^{-\gamma} \left ( t \vee t_1 \right )^{-1} \left ( h^{m-1} (t_1) + h^{m-1}(t_2) \right )$$
$$\leq 2 \int_1^{\infty} dt_1 \ t_1^{-\gamma} \left ( t \vee t_1 \right )^{-1} \int_{t_1}^{\infty}
dt_2 \ t_2^{-\gamma} \left ( t \vee t_1 \vee t_2 \right )^{-1} \ h^{m-1} (t_1) \leq 2 Q_m(t)
\eqno(3.50)$$

\noi since $h$ is decreasing in $t$ and
$$\int_{t_1}^{\infty} dt_2 \ t_2^{-\gamma} \left ( t \vee t_1 \vee t_2 \right )^{-1} \leq h\left ( t
\vee t_1 \right ) \leq h(t_1) \quad .$$

\noi Now (3.46) follows from (3.49) and (3.45). \par \vskip 3 truemm

(3.47). We first define for future use
$$Q_m = Q_m^- + Q_m^+ = t^{-1} \int_1^t dt_1 \ t_1^{-\gamma} \ h^m(t_1) + \int_t^{\infty} dt_1 \
t_1^{-1-\gamma} \ h^m (t_1) \quad , \eqno(3.51)$$
$$P_m = P_m^- + P_m^+ = h(t) \int_1^t dt_1 \ t_1^{-\gamma} \ h^m(t_1) + \int_t^{\infty} dt_1 \
t_1^{-\gamma} \ h^{m+1} (t_1) \quad , \eqno(3.52)$$
$$R_m = R_m^- + R_m^+ = \int_t^{\infty} dt_1 \ t_1^{-2} \ P_m^-(t_1) + \int_t^{\infty} dt_1 \
t_1^{-2} \ P_m^+(t_1) \eqno(3.53)$$

\noi and we estimate $R_m^-$ and $R_m^+$ separately. We first estimate
$$R_m^- = \int_t^{\infty} dt_1 \ t_1^{-2} \ h(t_1) \int_1^{t_1} dt_2 \ t_2^{-\gamma} \ h^m(t_2)$$
$$\leq h(t) \int_1^{\infty} dt_2 \ t_2^{-\gamma} \left ( t \vee t_2 \right )^{-1} \ h^m (t_2) =
h(t) \ Q_m(t) \eqno(3.54)$$
\noi by the monotony of $h$ and after performing the integral over $t_1$. We next use the
differential equation
$$\gamma h + t\ h' = t^{-1} \eqno(3.55)$$
\noi satisfied by $h$ to rewrite $P_m^+$ as follows
$$\gamma P_m^+(t) = - \int_t^{\infty} dt_1 \ t_1^{1 - \gamma} \ h^m(t_1) \ h'(t_1) +
\int_t^{\infty} dt_1 \ t_1^{-1-\gamma} \ h^m(t_1) \quad .$$

\noi Integrating by parts in the first integral and using (3.51), we obtain 
$$(m + 1)\gamma \ P_m^+(t) = t^{1 - \gamma} \ h^{m+1}(t) + (1 - \gamma) P_m^+(t) + (m + 1) \
Q_m^+(t)$$

\noi namely
$$\left ( (m + 2) \gamma - 1 \right ) P_m^+(t) = t^{1 - \gamma} \ h^{m+1}(t) + (m + 1) \ Q_m^+(t)
\quad .$$

\noi Substituting that result into the definition of $R_m^+(t)$, we obtain 
$$\left ( (m + 2) \gamma - 1 \right ) R_m^+(t) = \int_t^{\infty} dt_1 \ t_1^{-1-\gamma} \
h^{m+1}(t_1) + (m + 1) \int_t^{\infty} dt_1 \ t_1^{-1 - \gamma} \ h^m(t_1) \left ( t^{-1} -
t_1^{-1} \right )$$
$$\leq Q_{m+1}^+ (t) + (m + 1) t^{-1} \ Q_m^+(t)$$
$$\leq \left ( h(t) + (m + 1) t^{-1} \right ) Q_m^+(t)$$

\noi by the monotony of $h$, 
$$\leq \left ( 1 + (m + 1) \gamma \right ) \ h \ Q_m^+(t) \eqno(3.56)$$

\noi by (3.24). Collecting (3.54) and (3.56) yields (3.47). \par \nobreak
\hfill $\sq$ 

\section{Cauchy problem and preliminary asymptotics for the auxiliary system}
\hspace*{\parindent}In this section, we collect a number of results from I on the Cauchy problem
and on the asymptotic behaviour of solutions for the auxiliary system

$$\hskip 1.7 truecm\left \{ \begin{array}{ll} \partial_t w = i\left ( 2t^2 \right )^{-1} \ \Delta w +
\left ( 2t^2 \right )^{-1} \ \left ( 2 \nabla \varphi \cdot \nabla + (\Delta \varphi )
\right ) w &\hskip 3 truecm (2.11)\equiv(4.1) \\ & \\ \partial_t \varphi = \left (2 t^2 \right
)^{-1} \ |\nabla \varphi |^2 + t^{-\gamma} \ g_0(w, w) \quad .  &\hskip 3 truecm (2.12)\equiv (4.2) 
\end{array}\right .$$

Those results are immediate extensions of results contained in I. The main differences are that
(i) the results are stated here in terms of $\varphi$ whereas they are stated in I in terms of $s
= \nabla \varphi$, and (ii) here we use systematically the estimating functions of time $h_0$ and
$h$ introduced in Section 3, thereby covering the whole interval $0 < \gamma \leq 1$. The
proofs will be sketched briefly or omitted. \par

We first recall the results on the local Cauchy problem with finite initial time (see Proposition
I.4.1). \\

\noi {\bf Proposition 4.1.} {\it Let $(k, \ell )$ be an admissible pair. Let $t_0 > 0$. Then for any
$(w_0, \varphi_0) \in H^k \oplus Y^{\ell}$, there exist $T_{\pm}$ with $0 \leq T_- < t_0 < T_+ \leq
\infty$ such that~:} \par

{\it (1) The system (4.1) (4.2) has a unique solution $(w , \varphi) \in {\cal C}(I, H^k \oplus
Y^{\ell})$ with $(w, \varphi) (t_0) = (w_0, \varphi_0)$, where $I = (T_- , T_+)$. If $T_- > 0$ (resp.
$T_+ < \infty$), then $|w(t)|_k + |\varphi (t)|_{\ell} \to \infty$ when $t$ decreases to $T_-$ (resp.
increases to $T_+$).} \par

{\it (2) If $(w_0, \varphi_0) \in H^{k'} \oplus Y^{\ell '}$ for some admissible pair $(k', \ell ')$
with $k' \geq k$ and $\ell ' \geq \ell$, then $(w, \varphi) \in {\cal C}(I, H^{k'} \oplus Y^{\ell
'})$.} \par

{\it (3) For any compact subinterval $J \subset \subset I$, the map $(w_0, \varphi_0) \to (w,
\varphi)$ is continuous from $H^{k-1} \oplus Y^{\ell - 1}$ to $L^{\infty}(J, H^{k-1} \oplus Y^{\ell -
1})$ uniformly on the bounded sets of $H^k \oplus Y^{\ell}$, and is pointwise continuous from $H^k
\oplus Y^{\ell}$ to $L^{\infty}(J,H^k \oplus Y^{\ell})$.} \\

We next recall the results on the local Cauchy problem in a neighborhood of infinity in time (see
Proposition I.5.1). \\

\noi {\bf Proposition 4.2.} {\it Let $(k , \ell )$ be an admissible pair. Let $(w_0,
\widetilde{\varphi}_0) \in H^k \oplus Y^{\ell}$ and define $a = |w_0|_k$ and $b =
|\widetilde{\varphi}_0|_{\ell}$. Then there exists $T_0 < \infty$, depending on $a$, $b$, such that
for all $t_0 \geq T_0$, there exists $T \leq t_0$, depending on $a$, $b$ and $t_0$, such that the
system (4.1) (4.2) with initial data $w(t_0)=w_0$, $\varphi (t_0) = h_0(t_0)
\widetilde{\varphi}_0$ has a unique solution $(w, \varphi )$ in the interval $[T, \infty )$ such
that $(w, h_0^{-1} \varphi ) \in ({\cal C} \cap L^{\infty})([T, \infty ), H^k \otimes Y^{\ell})$.
One can define $T_0$ and $T$ by}
$$C  ( b + a^2 ) \ h \left ( T_0 \right  ) = 1 \eqno(4.3)$$
$$T = h_0 \left ( t_0 \right ) \ h \left ( T_0 \right )^{-1}  \eqno(4.4)$$

\noi {\it and the solution $(w, \varphi )$ is estimated for all $t \geq T$ by}
$$|w(t)|_k \leq C \ a \eqno(4.5)$$
$$|\varphi (t) |_{\ell} \leq C ( b + a^2 ) \ h_0 \left ( t \vee t_0 \right ) \quad
. \eqno(4.6)$$
\vskip 3 truemm

\noi {\bf Sketch of proof.} The proof is almost identical with that of Proposition I.5.1 and
follows from a priori estimates of the maximal solution obtained from Proposition 4.1. Define
$y = |w|_k$ and $z = |\varphi |_{\ell}$. By Lemmas 3.2 and 3.3, $y$ and $z$ satisfy 
$$ \left \{ \begin{array}{l} |\partial_t y | \leq C \ t^{-2} \ y \ z \\ \\ |\partial_t z| \leq C
\ t^{-2} \ z^2 + C \ t^{-\gamma} \ y^2 \quad . \end{array} \right . \eqno(4.7)$$

\noi For $t \geq t_0$, we take $\bar{t} > t_0$, we define $Y \equiv Y(\bar{t}) = \parallel y;
L^{\infty} ([t_0, \bar{t}\, ])\parallel$ and $Z \equiv Z(\bar{t}) =$\break \noindent $\parallel h_0
(t)^{-1} z; L^{\infty} ([t_0, \bar{t}])\parallel $, we substitute those definitions into (4.7), we
integrate over $t$ with the appropriate initial condition and we obtain
$$ \left \{ \begin{array}{l} Y \leq a + C \ Y \ Z \ h(t_0) \\ \\ Z \leq b + C \ Y^2 + C \ Z^2
\ h(t_0)  \end{array} \right . \eqno(4.8)$$

\noi by (3.19) (3.22). \par

For $t \leq t_0$, we take $\bar{t} < t_0$, we define $Y \equiv Y(\bar{t}) = \parallel y;L^{\infty}
([\bar{t},t_0])\parallel$ and $Z \equiv Z(\bar{t}) =$ \break \noindent $\parallel z;
L^{\infty}([\bar{t}, t_0])\parallel$, we substitute those definitions into (4.7), we integrate over
$t$ with the ap\-pro\-pria\-te initial condition and we obtain
$$ \left \{ \begin{array}{l} Y \leq a + C \ t^{-1} \ Y \ Z  \\ \\ Z \leq \left ( b + C \ Y^2
\right )\ h_0 (t_0) + C \ t^{-1} \ Z^2 \quad .  \end{array} \right . \eqno(4.9)$$

The proof then proceeds from (4.8) and (4.9) in the same way as that of Proposition I.5.1. \par
\nobreak
\hfill $\sq$

For subsequent applications, we shall need the following lemma, which is essentially identical with
Lemma I.5.1. \\

\noi {\bf Lemma 4.1.} {\it Let $a > 0$, $b > 0$, $t_0 > 1$ and let $y$, $z$ be nonnegative
continuous functions satisfying $y(t_0) = y_0$, $z(t_0) = z_0$ and}
$$ \left \{ \begin{array}{l} |\partial_t y | \leq t^{-2} \ h_0 (t) \ b \ y + t^{-2} \ a\ z \\ \\
|\partial_t z | \leq t^{-2} \ h_0(t) \ b \ z + t^{-\gamma} \ a\ y  \quad .  \end{array} \right .
\eqno(4.10)$$

\noi {\it Define $\bar{y}$, $\bar{z}$ by}
$$(y, z) = (\bar{y}, \bar{z}) \exp  ( b|h(t) - h(t_0)|) \quad . \eqno(4.11)$$

\noi {\it Then for $\gamma (t_0^{\gamma} \wedge t^{\gamma}) \geq 2a^2$, the following estimates
hold~:} $$ \left \{ \begin{array}{l} \bar{y} \leq 2 ( y_0 + a\ z_0 \ t_0^{-1} )   \\ \\
\bar{z} \leq z_0 + 2a  ( y_0 + a \ z_0 \ t_0^{-1} ) \ h_0(t)    \end{array}
\right . \eqno(4.12)$$
\noi {\it for $t \geq t_0$, and}
$$ \left \{ \begin{array}{l}  \bar{y} \leq y_0 + 2a \left ( z_0 + a\ y_0 \ h_0(t_0) \right ) \
t^{-1}  \\ \\  \bar{z} \leq 2 \left ( z_0 + a \ y_0 \ h_0 (t_0) \right )   \end{array} \right
. \eqno(4.13)$$

\noi {\it for $1 \leq t \leq t_0$.} \\

As an easy consequence of Lemma 4.1, we obtain the following uniqueness result at infinity for the
system (4.1) (4.2) (see Proposition I.5.2). \\

\noi {\bf Proposition 4.3.} {\it Let $(k, \ell )$ be an admissible pair. Let $(w_i, \varphi_i)$, $i
= 1,2$ be two solutions of the system (4.1) (4.2) such that $(w_i, h_0^{-1} \varphi_i ) \in
L^{\infty}([T, \infty ), H^k \oplus Y^{\ell})$ for some $T > 0$ and such that $|w_1(t) -
w_2(t)|_{k-1} \ h_0(t)$ and $|\varphi_1(t) - \varphi_2 (t)|_{\ell -1}$ tend to zero when $t \to
\infty$. Then $(w_1, \varphi_1) = (w_2, \varphi_2)$.}\\

We finally recall the existence result for the limit of $w(t)$ as $t \to \infty$ for the solutions
of the system (4.1) (4.2) obtained in Proposition 4.2 (see Proposition I.5.3). \\

\noi {\bf Proposition 4.4.} {\it Let $(k, \ell )$ satisfy $k \leq \ell + 1$ and $\ell > n/2$. Let
$(w, \varphi )$ satisfy (4.1) and be such that $(w, h_0^{-1} \varphi ) \in ({\cal C} \cap
L^{\infty}) ([T, \infty ), H^k \oplus Y^{\ell})$ for some $T >0$. Let}
$$a = \parallel w; L^{\infty}([T, \infty ), H^k)\parallel \quad , \quad b = \parallel h_0^{-1}
\varphi ; L^{\infty}([T, \infty ), Y^{\ell})\parallel \quad . \eqno(4.14)$$

\noi {\it Then there exists $w_+ \in H^k$ such that $w(t)$ tends to $w_+$ strongly in $H^{k-1}$
and weakly in $H^k$ when $t \to \infty$. Furthermore the following estimates hold}
$$|w_+|_k \leq a \eqno(4.15)$$
$$|w(t_0) - w(t)|_{k-1} \leq C \ a\ b \ h(t_0 \wedge t) \eqno(4.16)$$
$$|w(t) - w_+|_{k-1} \leq C \ a\ b \ h(t) \eqno(4.17)$$

\noi {\it for $t_0$, $t$ sufficiently large, namely $bh(t_0 \wedge t) \leq C$ or $bh(t) \leq C$.}

\section{Existence and properties of the asymptotic dynamics}
\hspace*{\parindent} In this section we derive the relevant properties of the solutions $(w_m,
\varphi_m)$ of the system (2.20) (2.21) with initial conditions (2.22) (2.23) and of the solutions
$(V, \chi )$ of the transport equations (2.25) (2.26) with initial conditions (2.27). We use
systematically the estimating functions of time $N_m$, $Q_m$ and $P_m$ defined by (3.25) (3.26)
and (3.31). We begin with the system (2.20) (2.21), which is solved by successive integrations, as
explained in Section 2. \\

\noi {\bf Proposition 5.1.} {\it Let $(k, \ell )$ be an admissible pair, let $p \geq 0$ be an
integer, let $w_+ \in H^{k+p}$ and let $a = |w_+|_{k+p}$. Let $\{w_0 = w_+, w_{m+1}\}$ and
$\{\varphi_m\}$, $0 \leq m \leq p$, be the solution of the system (2.20) (2.21) with initial
conditions (2.22) (2.23). Then} \par
{\it (1) $w_{m+1} \in {\cal C}([1, \infty ), H^{k+p-m-1})$, $\varphi_m \in {\cal C}([1, \infty ),
Y^{\ell + p - m})$ and the following estimates hold for all $t \geq 1$~:}
$$|w_{m+1}(t)|_{k+p-m-1} \leq A(a) \ Q_m(t) \eqno(5.1)$$
$$|\varphi_m(t)|_{\ell + p-m} \leq A(a) \ N_m(t) \eqno(5.2)$$
\noi {\it for some estimating function $A(a)$.} \par
{\it If in addition $(p + 2)\gamma > 1$ and if we define $\varphi_{p+1}$ by (2.21) with initial
condition $\varphi_{p+1}(\infty ) = 0$, then $\varphi_{p+1} \in {\cal C}([1, \infty ), Y^{\ell -
1})$ and the following estimate holds~:}
$$|\varphi_{p+1}(t)|_{\ell - 1} \leq A(a) \ P_p(t) \quad . \eqno(5.3)$$

{\it (2) The functions $\{\varphi_m \}$ are gauge invariant in the following sense. If $w'_+ =
w_+ \exp (i \sigma )$ for some real valued function $\sigma$ and if $w'_+$ gives rise to
$\{\varphi '_m\}$, then $\varphi '_m = \varphi_m$ for $0 \leq m \leq p+1$.} \par

{\it (3) The map $w_+ \to \{w_{m+1}, \varphi_m\}$ is uniformly Lipschitz continuous on the
bounded sets from the norm topology of $w_+$ in $H^{k+p}$ to the norms $\parallel Q_m^{-1} \
w_{m+1};L^{\infty}([1, \infty ), H^{k+p-m-1})\parallel$ and $\parallel N_m^{-1} \
\varphi_m;L^{\infty}([1, \infty ), Y^{\ell + p-m})\parallel$, $0 \leq m \leq p$. A similar
continuity holds for $\varphi_{p+1}$.} \\

\noi {\bf Proof.} {\bf Part (1).} The proof proceeds by induction on $m$. We assume the results to
hold for $(w_j, \varphi_j)$ for $j \leq m$ and we prove them for $w_{m+1}$ and $\varphi_{m+1}$.
We first consider $w_{m+1}$ which is obtained from (2.20). From Lemma 3.2, especially (3.6) with
$(k, \ell)$ replaced by $(k+p-m,\ell +p-m)$ which is again an admissible pair and from the
induction assumption, we obtain
$$\left | \partial_t \ w_{m+1} \right |_{k+p-m-1} \leq A(a) \ t^{-2} \Big \{ \sum_{0 \leq j
\leq m-1} N_j (t) \ Q_{m-j-1}(t) + N_m(t) \Big \} \quad . \eqno(5.4)$$

\noi Integrating (5.4) between $t$ and $\infty$, using the initial condition $w_{m+1}(\infty ) =
0$ and using (3.39) (3.35) shows that $w_{m+1} \in {\cal C}([1, \infty ), H^{k+p-m-1})$ and
that $w_{m+1}$ satisfies (5.1). \par

We next consider $\varphi_{m+1}$ which is obtained from (2.21). From Lemma 3.2, especially
(3.7) (3.9) with again $(k, \ell )$ replaced by $(k+p-m,\ell +p-m)$, from the induction
assumption and from the result for $w_{m+1}$, we obtain 
$$\left | \partial_t \ \varphi_{m+1} \right |_{\ell + p-m-1} \leq A(a) \Big \{ t^{-2} \sum_{0
\leq j \leq m} N_j(t) \ N_{m-j}(t) + t^{-\gamma} \Big ( \sum_{0\leq j \leq m - 1} Q_j(t) \
Q_{m-1-j}(t) + Q_m(t)\Big ) \Big \} \ . \eqno(5.5)$$

\noi Integrating (5.5) between 1 and $t$, using the initial condition $\varphi_{m+1}(1) = 0$,
and using (3.38) (3.36) and (3.40) (3.43) (3.44) shows that $\varphi_{m+1} \in {\cal C}([1,
\infty ), Y^{\ell + p-m-1})$ and that $\varphi_{m+1}$ satisfies (5.2). \par

We finally assume that $(p + 2)\gamma > 1$ and estimate $\partial_t \ \varphi_{p+1}$ by (5.5)
with $m = p$. The last result and in particular the estimate (5.3) then follow by integration
between $t$ and $\infty$ and use of (3.38) (3.37) and (3.40) (3.41) (3.42) with $m = p$. \\

\noi {\bf Part (2).} We define for $0 \leq m \leq p+1$
$$B_m = \sum_{0\leq j \leq m} \bar{w}_j \ w_{m-j}$$

\noi so that $B_0 = |w_+|^2$ and $B_m$ is bounded in time and tends to zero at infinity for $m \geq
1$, for instance in $H_1^1$ norm. The equation (2.21) for $\varphi_{m+1}$ can be rewritten as
$$\partial_t \ \varphi_{m+1} = ( 2t^2 )^{-1} \sum_{0 \leq j \leq m} \nabla \varphi_j
\cdot \nabla \varphi_{m-j} + t^{-\gamma} \lambda \ \omega^{\mu - n} \ B_{m+1} \quad . \eqno(5.6)$$

\noi We next compute

$$\begin{array}{ll}\partial_t \ B_{m+1} &=  ( 2t^2 )^{-1}  \displaystyle{\sum\limits_{0
\leq j \leq m}} 2{\rm Re} \ \bar{w}_j \displaystyle{\sum\limits_{0 \leq i \leq m - j}}  ( 2
\nabla \varphi_{m-i-j} \cdot \nabla +  ( \Delta \ \varphi_{m-i-j} ) ) w_i \\ & \\
&= t^{-2} \displaystyle{\sum\limits_{0 \leq k \leq m}}  ( \nabla \varphi_{m-k} \cdot \nabla +
 ( \Delta \ \varphi_{m-k} )  ) B_k \quad .\end{array} \eqno(5.7)$$

\noi Using (5.6) and (5.7), we now show by induction on $m$ that $B_m$ and $\varphi_m$ are gauge
invariant. In fact assume that $B_j$ and $\varphi_j$ are gauge invariant for $j \leq m$. Then
$\partial_t \ B_{m+1}$ is gauge invariant by (5.7) and therefore $B_{m+1}$ is gauge invariant
because $B_{m+1}(\infty ) = 0$. Substituting that result into (5.6) and using the induction
assumption, we obtain from (5.6) that $\partial_t \ \varphi_{m+1}$ is gauge invariant, and
therefore $\varphi_{m+1}$ is gauge invariant since $\varphi_{m+1} (1) = 0$ for $m < p$ and
$\varphi_{p+1}(\infty ) = 0$. \\

\noi {\bf Part (3).} Let $\{w_m, \varphi_m\}$ and $\{w'_m, \varphi '_m\}$ be the solutions of the
system (2.20) (2.21) associated with $w_+$ and $w'_+$. From the fact that the RHS of (2.20) (2.21)
are bilinear, it follows as in Part (1) by induction on $m$ that the following estimates hold,
with $a = |w_+|_{k+p} \vee |w'_+|_{k+p}$~: 
$$\left | w_{m+1} - w'_{m+1} \right |_{k+p-m-1} \leq A(a)\  |w_+ - w'_+|_{k+p} \ Q_m(t)
\eqno(5.8)$$
$$\left | \varphi_{m} - \varphi '_{m} \right |_{\ell +p-m} \leq A(a)\  |w_+ - w'_+|_{k+p} \ N_m(t)
\eqno(5.9)$$

\noi for $0 \leq m \leq p$, and if $(p + 2)\gamma > 1$,
$$\left | \varphi_{p+1} - \varphi'_{p+1} \right |_{\ell -1} \leq A(a)\  |w_+ - w'_+|_{k+p} \ P_p(t)
\quad . \eqno(5.10)$$

The continuity as stated in Part (3) follows from those estimates. \par \nobreak
\hfill $\sq$ 

\noi {\bf Remark 5.1.} There is no upper bound on $p$ in Proposition 5.1. However if $(p+1) \gamma
> 1$, all the $(w_m, \varphi_m)$ with $(m+1) \gamma > 1$ have the same asymptotic behaviour in time
and behave respectively as $t^{-1}$ and Constant as $t \to \infty$, because $Q_m$ and $N_m$
saturate to those behaviours in that case. \\

We define for future reference (see also (2.16) (2.17))
$$W_m = \sum_{0 \leq j \leq m} w_j \quad , \quad \phi_m = \sum_{0 \leq j \leq m} \varphi_j
\eqno(5.11)$$

\noi where $w_j$, $\varphi_j$ are obtained by Proposition 5.1. \par

We now turn to the study of the transport equation 
$$\partial_t \ V =  ( 2t^2 )^{-1} \ \left ( 2 \nabla \phi \cdot \nabla + (\Delta \phi )
\right ) V \eqno(5.12)$$

\noi which we shall use later with $\phi = \phi_{p-1}$, as explained in Section 2 (see (2.25)). \\

\noi {\bf Proposition 5.2.} {\it Let $\ell > n/2$ and $1 \leq k \leq \ell$. Let $T \geq 1$, $I =
[T, \infty )$, let $\phi \in {\cal C} (I, Y^{\ell})$ with $h_0^{-1}\phi \in L^{\infty} (I,
Y^{\ell})$ and let $w_+ \in H^k$. Then} \par
{\it (1) The equation (5.12) has a solution $V \in ({\cal C} \cap L^{\infty})(I,H^k)$ which is
estimated by}
$$\parallel V: L^{\infty} (I, H^k) \parallel \ \leq |w_+|_k \exp  ( C \ b \ \gamma^{-1} )
\eqno(5.13)$$ \noi {\it where}
$$b = \parallel h_0^{-1} \phi ; L^{\infty} (I, Y^{\ell})\parallel \quad , \eqno(5.14)$$

\noi {\it and which tends to $w_+$ at infinity in the sense that}
$$|V(t) - w_+|_{k-1} \leq C \ b \exp ( C \ b \ \gamma^{-1} ) \ |w_+|_k \ h(t) \quad .
\eqno(5.15)$$

{\it (2) The solution $V$ is unique in $L^{\infty} (I, L^2)$ under the condition that $\parallel
V(t) - w_+\parallel_2$ tends to zero as $t \to \infty$.} \par

{\it (3) The map $(w_+ , \phi ) \to V$ is uniformly Lipschitz continuous in $w_+$ for the norm
topology of $H^k$ and is continuous in $\phi$ for the topology of convergence in $Y^{\ell}$
pointwise in $t$ to the norm topology of $L^{\infty}(I, H^k)$ for $h_0^{-1} \phi$ in bounded sets
of $L^{\infty}(I, Y^{\ell})$.} \\

\noi {\bf Proof.} {\bf Part (1).} We first take $t_0 \in I$. Using a regularization (for instance
parabolic), energy estimates as in Lemmas 3.2 and 3.3 (see especially (3.6) and (3.12)), and a
limiting procedure, one obtains easily the existence of a solution $V_{t_0}$ of the equation (5.12)
with initial condition $V_{t_0}(t_0) = w_+$, and such that 
$$V_{t_0} \in {\cal C} ( I, H^{k-1} ) \cap ( {\cal C}_w \cap L^{\infty} ) ( I, H^k  ) \quad .$$

Using the same energy estimates, one then shows that
$$\Big | \partial_t | V_{t_0}(t)|_k \Big | \leq C \ b\ t^{-2} \ h_0(t) \left | V_{t_0}(t) \right
|_k \eqno(5.16)$$
$$\Big | \partial_t | V_{t_0}(t) - w_+|_{k-1} \Big | \leq C \ b\ t^{-2} \ h_0(t) \left ( \Big |
V_{t_0}(t) - w_+ \Big |_{k-1} + |w_+|_k \right ) \eqno(5.17)$$

\noi and for two solutions $V_{t_0}$ and $V_{t_1}$ associated with $t_0$ and $t_1$
 $$\Big | \partial_t | V_{t_0}(t) - V_{t_1}(t)|_{k-1} \Big | \leq C \ b\ t^{-2} \ h_0(t) 
\Big | V_{t_0}(t) - V_{t_1}(t) \Big |_{k-1} \quad . \eqno(5.18)$$

Integrating (5.16) (5.17) between $t_0$ and $t$ and integrating (5.18) between $t_1$ and $t$, we
obtain respectively
$$|V_{t_0}(t)|_k \leq |w_+|_k \exp \left ( C \ b |h(t) - h(t_0)| \right ) \leq |w_+|_k \exp (
C \ b \ \gamma^{-1} ) \quad , \eqno(5.19)$$
 $$\begin{array}{ll} \Big | V_{t_0}(t) - w_+ \Big |_{k-1} &\leq |w_+|_k \left ( \exp \left ( C \ b
|h(t) - h(t_0)| \right ) - 1 \right ) \\ &\\ &\leq |w_+|_k
\ C \ b \exp \left ( C \ b \ \gamma^{-1} \right ) |h(t) - h(t_0)| \quad ,   \end{array}
\eqno(5.20)$$      

$$\begin{array}{ll} \Big | V_{t_0}(t) - V_{t_1}(t) \Big |_{k-1} &\leq \Big | V_{t_0}(t_1) -
w_+ \Big |_{k-1} \exp \left ( C \ b |h(t) - h(t_1)| \right ) \\ & \\ &\leq \Big |
V_{t_0}(t_1) - w_+ \Big |_{k - 1} \exp \left ( C \ b \ \gamma^{-1} \right ) \quad .   \end{array}
\eqno(5.21)$$

\noi Substituting (5.20) into (5.21) yields
$$ \left | V_{t_0}(t) - V_{t_1}(t) \right |_{k-1} \leq | w_+  |_{k} \ C \ b \exp ( 2 C \ b
\gamma^{-1} ) \left | h(t_1) - h(t_0) \right | \quad . \eqno(5.22)$$

\noi It follows from (5.22) that when $t_0 \to \infty$, $V_{t_0}$ has a limit $V \in ({\cal C}
\cap L^{\infty})(I, H^{k-1})$ satisfying (5.15). One sees easily that $V$ satisfies the equation
(5.12). From the estimate (5.19) it follows by a standard compactness argument that $V \in
({\cal C}_w \cap L^{\infty})(I, H^k)$ and that $V$ satisfies the estimate (5.13). Furthermore
$V$ also satisfies (5.16) so that $|V(t)|_k$ is Lipschitz continuous in $t$, which together
with weak continuity in $H^k$ implies strong continuity in $H^k$. \\

\noi {\bf Part (2).} If $V_1$ and $V_2$ are two solutions of (5.12) one obtains by the same
energy estimates as above
$$\parallel V_1(t) - V_2(t)\parallel_2 \ \leq \ \parallel V_1(t') - V_2(t')\parallel_2 \exp
\left ( C \ b |h(t) - h(t')|\right ) \quad . \eqno(5.23)$$

\noi Taking the limit $t' \to \infty$ shows that $V_1 = V_2$. \\

\noi {\bf Part (3).} Continuity of $V$ with respect to $w_+$ follows immediately from the
linearity of the equation (5.12) and from the estimate (5.13). In order to prove continuity with
respect to $\phi$, we first derive an estimate for the difference of two solutions $V_1$ and
$V_2$ associated with $\phi_1$ and $\phi_2$. We assume that $\phi_1 \in {\cal C}(I, Y^{\ell +
1})$ with $h_0^{-1} \phi \in L^{\infty}(I, Y^{\ell + 1})$, that $\phi_2 \in {\cal C}(I,
Y^{\ell})$ with $h_0^{-1} \phi_2 \in L^{\infty}(I, Y^{\ell})$, that $V_1 \in ({\cal C} \cap
L^{\infty})(I, H^{k+1})$ and that $V_2 \in ({\cal C} \cap L^{\infty})(I, H^k)$. Let $V_- = V_1
- V_2$ and $\phi_- = \phi_1 - \phi_2$. It follows from (5.12) that
$$\partial_t V_- = ( 2t^2 )^{-1} \Big \{ \left ( 2 \nabla \phi_2 \cdot \nabla +
(\Delta \phi_2) \right ) V_- + \left ( 2 \nabla \phi_- \cdot \nabla + (\Delta \phi_-)\right )
V_1 \Big \} \quad . \eqno(5.24)$$

\noi Let
$$ a = \displaystyle{\mathrel{\mathop {\rm Max}_{i=1,2}}}\parallel V_i ; L^{\infty}(I, H^k)\parallel
\quad , \quad b = \displaystyle{\mathrel{\mathop {\rm Max}_{i=1,2}}}\parallel h_0^{-1} \phi_i ;
L^{\infty}(I, Y^{\ell}) \parallel \quad . $$

\noi Estimating (5.24) by Lemma 3.5, we obtain
$$\begin{array}{ll} \Big | \partial_t| V_-|_k \Big | &\leq C \ t^{-2} \Big \{
|\phi_2|_{\ell} \ |V_-|_k + |\phi_-|_{\ell} \ |V_1|_k + |\phi_-|_* \ |V_1|_{k+1} \Big \} \\
& \\ &\leq C \ t^{-2} \Big \{ b\ h_0 \ |V_-|_k + a |\phi_-|_{\ell} + |\phi_-|_* \ |V_1|_{k+1}
\Big \}  \end{array} \eqno(5.25)$$

\noi where $|f|_* \equiv \parallel \nabla f\parallel_{\infty}$. On the other hand, by Lemma
3.4 we obtain
$$\begin{array}{ll} \Big | \partial_t| V_1|_{k+1} \Big | &\leq C \ t^{-2} \Big \{
|\phi_1|_{\ell} \ |V_1|_{k+1} + |\phi_1|_{\ell +1} \ |V_1|_k \Big \} \\ &\\
&\leq C \ t^{-2} \Big \{ b\ h_0 \ |V_1|_{k+1} + a |\phi_1|_{\ell +1}  \Big \} \quad . \end{array}
\eqno(5.26)$$

\noi Integrating (5.26) between $t_0$ and $t$ and using the fact that 
$$|\partial_t y | \leq C_0 \ t^{-2} \ h_0\ y + z \eqno(5.27)$$

\noi implies
$$\begin{array}{ll} y(t) &\leq y(t_0) \exp \left ( C_0 \left | h(t) - h(t_0) \right | \right ) +
\left  | \displaystyle{\int_{t_0}^t} dt_1 \ z(t_1) \exp \left ( C_0 \left | h(t) - h(t_1)\right |
\right ) \right | \\ &\\ &\leq
\exp \left (  C_0 \ \gamma^{-1} \right ) \left ( y(t_0) + \left | \displaystyle{\int_{t_0}^t} dt_1 \
z(t_1) \right |  \right ) \quad , \end{array} \eqno(5.28)$$

\noi we obtain
$$\left | V_1(t) \right |_{k+1} \leq C \left ( \left | V_1 (t_0) \right |_{k+1} + \left |
\int_{t_0}^t dt_1 \ t_1^{-2} \ \left | \phi_1 (t_1) \right |_{\ell + 1} \right | \right ) 
\eqno(5.29)$$

\noi where $C$ depends on $a$, $b$. Substituting (5.29) into (5.25) and integrating between
$t_0$ and $t$ yields similarly
$$|V_-(t)|_k \leq C \left \{ |V_-(t_0)|_k + \left | \int_{t_0}^t dt_1 \ t_1^{-2}
|\phi_-(t_1)|_{\ell} \right | + | V_1(t_0)|_{k+1} \left | \int_{t_0}^t dt_1 \ t_1^{-2}
|\phi_-(t_1)|_* \right | \right .$$
$$\left . + \left | \int_{t_0}^t dt_1 \ t_1^{-2} \ |\phi_-(t_1)|_* \int_{t_0}^{t_1} dt_2 \ t_2^{-2} \
|\phi_1(t_2) |_{\ell + 1} \right | \right \} \quad . \eqno(5.30)$$

\noi In particular if $V_1$ and $V_2$ are solutions of (5.12) with $\phi_1$ and $\phi_2$
respectively and with initial data $w_{+1} \in H^{k+1}$ and $w_{+2} \in H^k$ at time $t_0$,
as obtained in Part (1), then the following estimate holds uniformly in $t_0$ and $t$
$$|V_-(t)|_k \leq C \left \{ \left | w_{+1} - w_{+2} \right |_k + \int_1^{\infty} dt \ t^{-2} \
|\phi_-(t) |_{\ell} \right .$$
$$\left . + \left ( |w_{+1}|_{k+1} + \int_1^{\infty} dt \ t^{-2} |\phi_1(t)|_{\ell + 1}
\right ) \int_1^{\infty} dt \ t^{-2} \ |\phi_-(t)|_* \right \} \ , \eqno(5.31)$$

\noi where all the integrals are convergent under the assumptions made on $\phi_1$ and
$\phi_2$. \par

We can now prove the continuity with respect to $\phi$. The proof proceeds as in Step 7 of
that of Proposition I.4.1. We introduce a regularization defined as follows. We choose a
function $\psi_1 \in {\cal S} ({I \hskip - 1 truemm R}^n)$ such that $\int dx \ \psi_1 (x) =
1$ and such that $|\xi|^{-2} (\widehat{\psi}_1 (\xi ) - 1)|_{\xi = 0} = 0$. We define
$\psi_{\varepsilon} (x) = \varepsilon^{-n} \psi_1(x/\varepsilon )$, so that
$\widehat{\psi}_{\varepsilon} (\xi ) = \widehat{\psi}_1 (\varepsilon \xi )$ and we define
the regularization by $f \to f_{\varepsilon} = \psi_{\varepsilon} * f$ for all $f \in
{\cal S}'$. An immediate computation yields
$$\parallel \partial f_{\varepsilon} \parallel_2 \ \leq \ \parallel \partial
\psi_{\varepsilon}\parallel_1 \ \parallel f \parallel_2 \ = \varepsilon^{-1} \parallel
\partial \psi_1 \parallel_1 \ \parallel f \parallel_2 \eqno(5.32)$$

\noi and
$$\parallel f_{\varepsilon} - f \parallel_{\infty} \ \leq \ \parallel (
\widehat{\psi}_{\varepsilon} - 1) \widehat{f} \parallel_1 \ \leq \varepsilon^{\theta} \parallel
|\xi|^{-n/2- \theta}  ( \widehat{\psi}_1 (\xi ) - 1 )  \parallel_2 \ \parallel f
;\dot{H}^{n/2+ \theta} \parallel \quad .  \eqno(5.33)$$

\noi Let now $w_+ \in H^k$ and $\phi$, $\phi ' \in {\cal C}(I, Y^{\ell})$ with $h_0^{-1} \phi$,
$h_0^{-1} \phi ' \in L^{\infty} (I, Y^{\ell})$ and such that 
$$\parallel h_0^{-1}\phi ; L^{\infty}(I, Y^{\ell}) \parallel \ \vee \ \parallel h_0^{-1} \phi ' ;
L^{\infty}(I, Y^{\ell}) \parallel \ \leq b \quad .$$

\noi Let $V$ and $V'$ be the solutions of the equation (5.12) with $\phi$ and $\phi '$ respectively
and with initial data $w_+$ at $t_0$ obtained in Part (1). We regularize $w_+$, $\phi$, $\phi '$ to
$w_{+ \varepsilon}$, $\phi_{\varepsilon}$, $\phi '_{\varepsilon}$, so that the following estimates
hold~:
$$|w_{+\varepsilon}|_{k+1} \leq C \ \varepsilon^{-1} |w_+|_k \quad , \quad
|\phi_{\varepsilon}|_{\ell + 1} \leq C \ \varepsilon^{-1} |\phi |_{\ell} \quad , \quad |\phi
'_{\varepsilon}|_{\ell + 1} \leq C \ \varepsilon^{-1} \ |\phi '|_{\ell} \eqno(5.34)$$

\noi and
$$|\phi - \phi_{\varepsilon}|_* \leq C \ \varepsilon^{3/2} \ |\phi |_{\ell} \quad , \quad |\phi ' -
\phi '_{\varepsilon}|_* \leq C \ \varepsilon^{3/2} \ |\phi '|_{\ell} \quad . \eqno(5.35).$$

\noi The estimates (5.34) follow from (5.32) and the estimates (5.35) follow from (5.33) and from
the definition of $Y^{\ell}$. Let $V_{\varepsilon}$ and $V'_{\varepsilon}$ be the solutions of
(5.12) obtained from $(w_{+ \varepsilon}, \phi_{\varepsilon})$ and $(w'_{+\varepsilon},
\phi'_{\varepsilon})$. We estimate
$$|V(t) - V'(t)|_k \leq |V(t) - V_{\varepsilon}(t)|_k + |V_{\varepsilon}(t) -
V'_{\varepsilon}(t)|_k + |V'_{\varepsilon} (t) - V'(t)|_k \quad . \eqno(5.36)$$

\noi We estimate the three norms in the RHS of (5.36) by applying successively (5.31) with $(V_1,
V_2) = (V_{\varepsilon}, V), (V_{\varepsilon}, V'_{\varepsilon})$ and $(V'_{\varepsilon}, V')$.
We obtain
$$|V(t) - V'(t)|_k \leq C \left \{ \left | w_+ - w_{+\varepsilon} \right |_k + \int_1^{\infty} dt
\ t^{-2} \Big ( \left | \phi(t) - \phi_{\varepsilon}(t) \right |_{\ell} + \left | \phi_{\varepsilon}(t)
- \phi '_{\varepsilon }(t) \right |_{\ell} \right .$$ $$\left . + \left |\phi '_{\varepsilon}(t) -
\phi '(t)\right |_{\ell} \Big ) + \varepsilon^{-1} \int_1^{\infty} dt \ t^{-2} \Big ( \left |
\phi(t) - \phi_{\varepsilon}(t) \right |_* + \left | \phi_{\varepsilon}(t) - \phi
'_{\varepsilon}(t)\right |_* + \left | \phi '_{\varepsilon}(t) - \phi '(t) \right |_*  \Big )
\right \} \eqno(5.37)$$

\noi where we have used (5.34). Using the inequalities
\begin{eqnarray*}
&&\left | \phi_{\varepsilon} - \phi '_{\varepsilon} \right |_{\ell} \leq \left | \phi
- \phi ' \right |_{\ell} \\
&&\left | \phi_{\varepsilon} - \phi '_{\varepsilon} \right |_* \leq \left | \phi - \phi ' \right
|_* \leq | \phi - \phi ' |_{\ell} \\
&&\left | \phi ' _{\varepsilon} - \phi ' \right |_{\ell} \leq \left | \phi_{\varepsilon} - \phi
\right |_{\ell} + 2 |\phi - \phi '|_{\ell} 
\end{eqnarray*}

\noi and (5.35), we can continue (5.37) by

$$\begin{array}{ll} |V(t) - V'(t)|_k &\leq C \Big \{ \left | w_+ - w_{+ \varepsilon} \right |_k +
\displaystyle{\int_1^{\infty}} dt \ t^{-2} \left | \phi (t) - \phi_{\varepsilon} (t) \right |_{\ell}
+ \varepsilon^{1/2} \\ &\\ &+ \left ( 1 + \varepsilon^{-1} \right ) \displaystyle{\int_1^{\infty}} dt
\ t^{-2} \left | \phi (t) - \phi '(t) \right |_{\ell} \Big \} \quad . \end{array} \eqno(5.38)$$ 
  
\noi For fixed $\phi$, by the Lebesgue dominated convergence theorem, the first integral in the RHS
tends to zero when $\varepsilon \to 0$, while the second integral tends to zero when $\phi ' \to
\phi$ in $Y^{\ell}$ pointwise in $t$. The RHS of (5.38) can then be made arbitrarily small by first
taking $\varepsilon$ sufficiently small and then letting $\phi '$ tend to $\phi$ in the previous
sense for fixed $\varepsilon$. \par \nobreak
\hfill $\sq$

We next turn to the analogous transport equation
$$\partial_t \chi = t^{-2} \ \nabla \phi \cdot \nabla \chi \eqno(5.39)$$

\noi which we shall use together with (5.12), as explained in Section 2. \\

\noi {\bf Proposition 5.3.} {\it Let $\ell > n/2$. Let $T \geq 1$, $I = [T, \infty )$, let $\phi
\in {\cal C}(I, Y^{\ell + 1})$ with $h_0^{-1} \phi \in L^{\infty}(I, Y^{\ell + 1})$ and let $\psi_+
\in Y^{\ell}$. Then} \par
{\it (1) The equation (5.39) has a solution $\chi \in ({\cal C} \cap L^{\infty})(I, Y^{\ell})$
which is estimated by}
$$\parallel \chi ; L^{\infty}(I, Y^{\ell})\parallel \ \leq |\psi_+|_{\ell} \exp ( C \ b \
\gamma^{-1} ) \eqno(5.40)$$

\noi {\it where}
$$b = \parallel h_0^{-1} \phi ; L^{\infty}(I, Y^{\ell + 1} ) \parallel \quad ,
\eqno(5.41)$$

\noi {\it and tends to $\psi_+$ at infinity in the sense that}
$$|\chi (t) - \psi_+ |_{\ell - 1} \leq C \ b \ \exp ( C \ b\ \gamma^{-1} )
|\psi_+|_{\ell} \ h(t) \quad . \eqno(5.42)$$

{\it (2) The solution $\chi$ is unique in $L^{\infty}(I,L^{\infty})$ under the condition that
$\parallel \chi (t) - \psi_+ \parallel_{\infty}$ tends to zero as $t \to \infty$.} \par

{\it (3) The map $(\psi_+ , \phi ) \to \chi$ is uniformly Lipschitz continuous in $\psi_+$ for the
norm topology of $Y^{\ell}$ and is continuous in $\phi$ for the topology of convergence in
$Y^{\ell + 1}$ pointwise in $t$ to the norm topology of $L^{\infty}(I, Y^{\ell})$, for $h_0^{-1}
\phi$ in bounded sets of $L^{\infty} (I, Y^{\ell + 1})$.} \par

{\it (4) Let in addition $w_+ \in H^k$ and let $V$ be the solution of (5.12) obtained in
Proposition 5.2. Then for fixed $\phi$, $V \exp (- i \chi )$ is gauge invariant in the following
sense~: if $(V, \chi )$ and $(V', \chi ')$ are the solutions obtained from $(w_+, \psi_+)$ and
$(w'_+, \psi '_+)$ and if $w_+ \exp (- i \psi_+) = w'_+ \exp (- i \psi '_+)$, then $V(t) \exp (- i
\chi (t)) = V'(t) \exp (-i\chi '(t))$ for all $t \in I$.} \\

\noi {\bf Proof.} {\bf Parts (1) (2) (3).} The proof is the same as that of Proposition 5.2,
starting from the estimates (3.13) (3.14) (3.16) (3.18) of Lemmas 3.3, 3.4 and 3.5. \\

\noi {\bf Part (4).} It follows from (5.12) and (5.39) that $V \exp (- i \chi )$ also satisfies
(5.12), with gauge invariant initial condition $V(\infty ) \exp (- i \chi ( \infty )) = w_+ \exp (-i
\psi_+)$. The result then follows from the uniqueness statement of Proposition 5.2, part (2). \par
\nobreak
\hfill $\sq$

In the subsequent applications, we shall use the solutions of the equations (5.12) and (5.39)
associated with $\phi = \phi_{p-1}$ defined by (5.11) (see (2.25) (2.26)). In particular we shall
use $V$ as a substitute for $W_p$, also defined by (5.11), and we shall need the fact that $V$ is a
sufficiently good approximation of $W_p$. We collect the relevant properties in the following
proposition. \\

\noi {\bf Proposition 5.4.} {\it Let $(k, \ell)$ be an admissible pair. Let $p \geq 1$ be an
integer. Let $w_+ \in H^{k+p+1}$, let $a = |w_+|_{k+p+1}$ and let $\phi = \phi_{p-1}$ be defined by
(5.11) and Proposition 5.1, so that $h_0^{-1} \phi \in ({\cal C} \cap  L^{\infty}) ([1, \infty
),Y^{\ell + 2})$. Let $V$ be the solution of (5.12) defined by Proposition 5.2, so that $V \in
({\cal C} \cap L^{\infty}) ([1, \infty ), H^{k+2})$.} \par

{\it (1) Let $W_p$ be defined by (5.11) and Proposition 5.1 so that $W_p \in ({\cal C} \cap
L^{\infty}) ([1, \infty ), H^{k+1})$. Then}
$$\left | V(t) - W_p(t) \right |_k \leq A(a) \ Q_p(t) \eqno(5.43)$$

\noi {\it for some estimating function $A(a)$.} \par

{\it (2) Let $\psi_+ \in Y^{\ell + 1}$ and let $\chi$ be the solution of (5.39) defined by
Proposition 5.3, so that $\chi \in ({\cal C} \cap L^{\infty})([1, \infty ), Y^{\ell + 1})$. Then}
$$\left | \chi (t) - \psi_+ \right |_{\ell} \leq A(a) \ |\psi_+|_{\ell + 1} \ h(t) \eqno(5.44)$$

\noi {\it and $V(t) \exp (-i \chi (t))$ is gauge invariant.}\par

{\it (3) The map $w_+ \to V$ is continuous from $H^{k+p+1}$ to $L^{\infty}([1, \infty ),
H^{k+2})$ and the map $(w_+, \psi_+) \to \chi$ is continuous from $H^{k+p+1} \oplus Y^{\ell + 1}$
to $L^{\infty}([1, \infty ), Y^{\ell + 1})$.}  \\

\noi {\bf Proof.} {\bf Part (1).} From (2.20) and (5.11) it follows that
$$\partial_t \ W_p = \left ( 2t^2 \right )^{-1} \Big \{ \left ( 2 \nabla \phi_{p-1} \cdot \nabla +
(\Delta \phi_{p-1}) \right ) W_p - \sum_{i \leq p - 1, j \leq p \atop{i+j \geq p}} \left ( 2 \nabla
\varphi_i \cdot \nabla + (\Delta \varphi_i) \right ) w_j \Big \}  \eqno(5.45)$$

\noi so that
$$\partial_t (V -  W_p) = \left ( 2t^2 \right )^{-1} \Big \{ \left ( 2 \nabla \phi_{p-1} \cdot
\nabla + (\Delta \phi_{p-1}) \right ) (V - W_p) + \sum_{i \leq p - 1, j \leq p \atop{i+j \geq p}}
\left ( 2 \nabla \varphi_i \cdot \nabla + (\Delta \varphi_i) \right ) w_j \Big \}\ . $$

\noi From Lemma 3.3, esp. (3.12) and Lemma 3.2, esp. (3.6), we obtain 
$$\Big | \partial_t| V - W_p |_k \Big | \leq C \Big \{ t^{-2} \ h_0(t) \ b \ |V - W_p|_k +
t^{-2} \sum_{i \leq p - 1, j \leq p \atop{i+j \geq p}} |\varphi_i|_{\ell} \ |w_j|_{k+1} \Big \}
\eqno(5.46)$$

\noi where
$$b = \parallel h_0^{-1} \ \phi_{p-1} ; L^{\infty}([1, \infty ), Y^{\ell})\parallel \quad .$$

\noi Integrating (5.46) between $t$ and infinity, using (5.1) (5.2) (5.27) (5.28), we obtain
$$\begin{array}{ll} \left |V(t) - W_p(t) \right |_k &\leq A(a) \displaystyle{\int_t^{\infty}} dt_1 \
t_1^{-2} \sum\limits_{i \leq p - 1, j \leq p \atop{i+j \geq p}} N_i(t_1) \ Q_{j-1}(t_1) \\ &\\
&\leq A(a) \sum\limits_{p\leq m \leq 2 p-1} Q_m(t) \leq A(a) \ Q_p(t) \end{array} \eqno(5.47)$$

\noi by (3.39) (3.35) and (3.28).\\

\noi {\bf Part (2)} is a partial rewriting of Proposition 5.3 in the special case $\phi =
\phi_{p-1}$. \\

\noi {\bf Part (3).} The continuity properties stated there follow by combining those of
Propositions 5.1, part (3), 5.2 part (3) and 5.3, part (3).  \par \nobreak \hfill $\sq$ 

\noi {\bf Remark 5.2.} By keeping track of the orders of derivation more accurately, one sees
easily that Proposition 5.4 holds with $(k, \ell )$ replaced everywhere by $(k - 1, \ell - 1)$. We
have stated Proposition 5.4 at the level of regularity which will be used in the subsequent
applications.        

\section{Asymptotics and wave operators for the auxiliary system}\nobreak
\hspace*{\parindent} In this section we derive the main technical results of this paper. We prove
that sufficiently regular solutions $(w, \varphi )$ of the auxiliary system (4.1) (4.2) have
asymptotic states $(w_+, \psi_+)$, and conversely that sufficiently regular asymptotic states $(w_+,
\psi_+)$ generate solutions $(w, \varphi )$ of the auxiliary system in the sense described in
Section 2, thereby allowing for the definition of the wave operator $\Omega_0 : (w_+, \psi_+) \to
(w, \varphi )$. \par

We first prove the existence of asymptotic states of solutions $(w, \varphi )$ of (4.1) (4.2). The
existence of $w_+$ is already established in Proposition 4.4 under rather general assumptions.
However the existence of $\psi_+$ requires a more complicated construction and stronger
assumptions. \\

\noi {\bf Proposition 6.1.} {\it Let $(k , \ell )$ be an admissible pair. Let $p \geq 0$ be an
integer. Let $T \geq 2$, $I = [T, \infty )$, let $(w, \varphi )$ be a solution of the system
(4.1) (4.2) such that $(w, h_0^{-1} \varphi ) \in ({\cal C} \cap L^{\infty})(I, H^{k+(p+1)\vee
2} \oplus Y^{\ell + p})$ and let}
$$a = \ \parallel w; L^{\infty}( [T, \infty ), H^{k+(p+1)\vee 2} ) \parallel \quad ,
\quad b = \ \parallel h_0^{-1}\varphi ; L^{\infty}( [T, \infty ), Y^{\ell + p} )
\parallel \quad .  \eqno(6.1)$$

\noi {\it Let $w_+ = \lim\limits_{t \to \infty} w(t) \in H^{k+p+1}$ be defined by Proposition 4.4.
Let $\{ w_{m+1} , \varphi_m\}$, $0 \leq m \leq p$ be defined by Proposition 5.1, and let $W_m$,
$\phi_m$, $0 \leq m \leq p$, be defined by (5.11). Then the following estimates hold for all $t \in
I$~:}  
$$\Big | w(t) - W_m(t) \Big |_{k+p-m-1} \ \leq A(a, b) \ Q_m(t) \eqno(6.2)$$
$$\Big | \varphi (t) - \phi_m(t) \Big |_{\ell + p - m -1} \ \leq A(a,b) \ N_{m+1}(t)
\eqno(6.3)$$
\noi {\it for $0 \leq m \leq p$, and for some estimating function $A(a,b)$.} \par
{\it If in addition $(p+2)\gamma > 1$, then the following limit exists}
$$\lim_{t \to \infty} \Big ( \varphi (t) - \phi_p(t) \Big ) = \psi_+ \eqno(6.4)$$

\noi {\it as a strong limit in $Y^{\ell - 1}$, and the following estimate holds}
$$\Big | \varphi (t) - \phi_p(t) - \psi_+ \Big |_{\ell - 1} \ \leq A(a, b) \ P_p(t) \quad . 
\eqno(6.5)$$   
\vskip 3 truemm
\noi {\bf Proof.} The proof proceeds by induction on $m$. For $0 \leq m \leq p$, we define 
$$q_{m+1} (t) = w(t) - W_m(t) \quad , \eqno(6.6)$$
$$\psi_{m+1}(t) = \varphi (t) - \phi_m(t) \quad . \eqno(6.7)$$

\noi We also define $q_0 = w$ and $\psi_0 = \varphi$. We assume that the estimates (6.2)
(6.3) hold for $(q_j , \psi_j)$, $0 \leq j \leq m$ and we derive them for $(q_{m+1} ,
\psi_{m+1})$. \par

We substitute the decompositions $w = W_m + q_{m+1}$ and $\varphi = \phi_m + \psi_{m+1}$ in
the LHS of (4.1) (4.2) and we partly substitute the decompositions $w = W_{m-1} + q_m$ and
$\varphi = \phi_{m-1} + \psi_m$ in the RHS of the same equations. Using in addition (2.20)
(2.21), we obtain 
$$\partial_t \ q_{m+1} =  ( 2t^2 )^{-1} \Big \{ i \Delta w + \left ( 2 \nabla
\varphi \cdot \nabla + (\Delta \varphi ) \right ) q_m$$
$$+ \left ( 2 \nabla \psi_m \cdot \nabla + (\Delta \psi_m )\right ) W_{m-1} + \sum_{0 \leq i,j \leq
m - 1\atop{i+j \geq m}} \left ( 2 \nabla \varphi_i \cdot \nabla + (\Delta \varphi_i) \right )
w_j \Big \} \quad ,  \eqno(6.8)$$        
$$\partial_t \ \psi_{m+1} = ( 2t^2 )^{-1} \Big \{ \left ( \nabla \varphi + \nabla
\phi_{m-1} \right ) \cdot \nabla \psi_m + \sum_{0 \leq i,j \leq m - 1\atop{i+j \geq m}} \nabla
\varphi_i \cdot \nabla \varphi_j \Big \}$$ 
$$ + t^{- \gamma} \Big \{ g_0(q_m, q_1) + g_0 (q_m, W_{m-1} - w_0) + 2 g_0 (q_{m+1}, w_0) + 
\sum_{0 \leq i,j \leq m - 1\atop{i+j \geq m+1}} g_0 (w_i, w_j) \Big \} \ . \eqno(6.9)$$

\noi The equation (6.9) holds only for $m \geq 1$ and the last bracket thereof has been obtained
by using the fact that
$$g_0\left ( q_m, w + W_{m-1} \right ) - 2g_0 \left ( w_m, w_0 \right ) = g_0(q_m, q_1) + g_0\left
( q_m, W_{m-1} - w_0 \right ) + 2g_0\left ( q_{m+1}, w_0 \right ) \quad .$$ 

\noi For $m = 0$, (6.9) should be replaced by
$$\partial_t \ \psi_1 = ( 2t^2 )^{-1} |\nabla \varphi |^2 + t^{-\gamma} ( g_0
(q_1, q_1) + 2g_0 (q_1, w_0) ) \quad . \eqno(6.9)_0$$

\noi We estimate the RHS of (6.8) (6.9) by Lemma 3.2 with $(k, \ell )$ replaced by $(k + p - m,
\ell + p - m)$, which is again an admissible pair. We use Proposition 5.1 to estimate $W_{m-1}$,
$w_j$, $\phi_{m-1}$ and $\varphi_j$, and we use the induction hypothesis to estimate $q_m$,
$q_1$, $q_{m+1}$ and $\psi_m$. Note that $W_{m-1} - w_0$, which occurs only for $m \geq 2$,
satisfies the same estimate as $q_1$. In the induction procedure, as in the proof of Proposition
5.1, one has first to complete the estimation of $q_{m+1}$ before estimating $\psi_{m+1}$. One
then obtains
$$\Big | \partial_t \ q_{m+1} \Big |_{k+p-m-1} \leq A(a, b) \ t^{-2} \Big \{ 1 + h_0 \ Q_{m-1} +
N_m + \sum_{m \leq i+j \leq 2(m-1)} N_i \ Q_{j-1} \Big \} \quad , \eqno(6.10)$$    
$$\Big | \partial_t \ \psi_{m+1} \Big |_{\ell +p-m-1} \leq A(a, b) \Big \{ t^{-2} \Big ( h_0
\ N_m + \sum_{m \leq i+j \leq 2(m-1)} N_i \ N_{j} \Big ) $$
$$ + t^{-\gamma} \Big ( h \ Q_{m-1} +
Q_m + \sum_{m+1 \leq i+j \leq 2(m-1)} Q_{i-1} \ Q_{j-1} \Big ) \Big \}\eqno(6.11)$$

\noi for $m \geq 1$, and
$$\hskip - 1.3 truecm \Big | \partial_t \ q_1 \Big |_{k+p-1} \leq A(a, b) \ t^{-2}(1 + h_0)
\eqno(6.10)_0$$     $$\Big | \partial_t \ \psi_1 \Big |_{\ell+p-1} \leq A(a, b) \left \{ t^{-2}h_0^2
+ t^{-\gamma} h \right \} \quad . \eqno(6.11)_0$$

\noi Integrating (6.10) between $t$ and infinity with the condition $q_{m+1}(\infty ) = 0$ which
follows from the definition and using (3.39) (3.35) (3.28) and the fact that $t^{-1} Q_m(1) \leq
Q_m(t)$ yields (6.2). \par

Integrating (6.11) between $T$ and $t$, using (3.38) (3.36) (3.27) for the first bracket and
(3.40) (3.43) (3.44) (3.28) for the second bracket yields 
\begin{eqnarray*}
\Big | \psi_{m+1}(t) \Big |_{\ell + p - m-1} &\leq& C + A(a,b) \ N_{m+1}(t) \\
&\leq& \left ( C \ N_{m+1}(T)^{-1} + A(a,b) \right ) \ N_{m+1}(t)
\end{eqnarray*}

\noi with
$$C = \Big | \psi_{m+1}(T) \Big |_{\ell + p - m -1} \quad , $$

\noi where we have used the fact that $N_{m+1}$ is increasing in $T$, and we have assumed that $T$ is
bounded away from 1. This yields (6.3). \par

We next turn to the proof of (6.5). In that case, the RHS of (6.11) with $m = p$ is integrable
in time, which proves the existence of the limit (6.4). Integrating (6.11) between $t$ and
infinity and using (3.38) (3.37) for the first bracket and (3.40) (3.41) (3.42) (3.28) for the
second bracket yields (6.5). \par \nobreak
\hfill $\sq$ \par

We now turn to the construction of solutions $(w, \varphi )$ of the system (4.1) (4.2) with
given asymptotic states $(w_+, \psi_+)$. For that purpose we first take a (large) positive $t_0$
and we construct the solution $(w_{t_0}, \varphi_{t_0})$ of (4.1) (4.2) with initial data
$(V(t_0), \phi_p(t_0) + \chi (t_0))$ at $t_0$. The solution $(w, \varphi )$ will then be
obtained therefrom by taking the limit $t_0 \to \infty$, as explained in Section 2. \\

\noi {\bf Proposition 6.2.} {\it Let $(k, \ell )$ be an admissible pair and let $p$ be an integer
such that $(p + 2)\gamma > 1$. Let $w_+ \in H^{k+(p+1)\vee 2}$ and $\psi_+ \in Y^{\ell + 1}$. Let
$\phi = \phi_{p-1}$ be defined by (5.11) and Proposition 5.1, so that $h_0^{-1} \phi \in ({\cal C}
\cap L^{\infty})([1, \infty ), Y^{\ell + 2})$. Let $V$ and $\chi$ be the solutions of (5.12) and
(5.39) respectively, obtained in Propositions 5.4, so that $(V, \chi ) \in ({\cal C} \cap
L^{\infty})([1, \infty ), H^{k+2} \oplus Y^{\ell + 1})$. Let}  
$$a_+ =  | w_+ |_{k+(p+1)\vee 2} \qquad , \qquad b_+ = |\psi_+|_{\ell + 1} \quad . \eqno(6.12)$$

\noi {\it Then there exist $T_0$ and $T$, $1 \leq T_0$, $T < \infty$, depending only on $(\gamma ,
p, a_+, b_+)$ such that for all $t_0 \geq T_0 \vee T$, the system (4.1) (4.2) with initial data
$w_{t_0} (t_0) = V(t_0)$, $\varphi_{t_0} (t_0) = \phi_p (t_0) + \chi (t_0)$ has a unique solution
in the interval $[T, \infty )$ such that $(w_{t_0}, h_0^{-1} \varphi_{t_0}) \in ({\cal C} \cap
L^{\infty}) ([T, \infty ), H^k \oplus Y^{\ell})$. One can define $T_0$ and $T$ by conditions of the
type}
$$A\left ( a_+, b_+ \right ) \ h(T_0) = 1 \eqno(6.13)$$
$$A\left ( a_+, b_+ \right ) \left ( (p + 2)\gamma - 1 \right )^{-1} \ h(T) = 1 \quad .
\eqno(6.14)$$

\noi {\it The solution satisfies the estimates} 
$$\Big | w_{t_0}(t) - V(t) \Big |_k \vee \Big |w_{t_0}(t) - W_p(t) \Big |_k \leq A \left ( a_+, b_+ \right
)\ Q_p(t_0) \eqno(6.15)$$
$$\Big | \varphi_{t_0}(t) - \phi_p(t) - \chi (t) \Big |_{\ell} \vee \Big |\varphi_{t_0}(t) -
\phi_p(t) - \psi_+\Big |_{\ell} \leq A \left ( a_+, b_+ \right )\ Q_p(t_0) \ h_0(t) \eqno(6.16)$$

\noi {\it for $t \geq t_0$,}
$$\Big | w_{t_0}(t) - V(t) \Big |_k \vee \Big |w_{t_0}(t) - W_p(t) \Big |_k \leq A \left ( a_+, b_+ \right
)\ Q_p(t) \eqno(6.17)$$
$$\Big | \varphi_{t_0}(t) - \phi_p(t) - \chi (t) \Big |_{\ell} \vee \Big |\varphi_{t_0}(t) -
\phi_p(t) - \psi_+\Big |_{\ell} \leq A \left ( a_+, b_+ \right )\ P_p(t) \eqno(6.18)$$

\noi {\it for $T \leq t \leq t_0$, and} 
$$|w_{t_0}(t)|_k  \leq A \left ( a_+, b_+ \right )  \quad , \quad |\varphi_{t_0}(t)|_{\ell}  \leq A
\left ( a_+, b_+ \right ) \ h_0(t) \eqno(6.19)$$

\noi {\it for all $t \geq T$.} \\

\noi {\bf Proof.} The result follows from Proposition 4.1 and standard globalisation arguments
provided we can derive (6.15) (6.16) (6.17) (6.18) as a priori estimates under the assumptions of
the proposition. Let $(w_{t_0}, \varphi_{t_0})$ be the maximal solution of (4.1) (4.2) with the
appropriate initial condition at $t_0$. Define $\widetilde{q} = w_{t_0} - V$ and
$\widetilde{\psi} = \varphi_{t_0} - \phi_p - \chi$. Comparing the equations (4.1) (4.2) and
(5.12) (5.39), we obtain 
$$\partial_t \ \widetilde{q} = ( 2t^2 )^{-1} \Big \{ i \Delta w_{t_0} + \left ( 2
\nabla \varphi_{t_0} \cdot \nabla + (\Delta \varphi_{t_0}) \right ) \widetilde{q} + \left ( 2
\nabla (\widetilde{\psi} + \varphi_p + \chi ) \cdot \nabla + (\Delta (\widetilde{\psi} +
\varphi_p + \chi ) )\right ) V \Big \} \eqno(6.20)$$
$$\partial_t \ \widetilde{\psi} = ( 2t^2 )^{-1} \Big \{ |\nabla \widetilde{\psi}|^2
+ 2 \nabla \widetilde{\psi}\cdot \nabla (\phi_p + \chi ) + |\nabla \chi |^2 + 2 \nabla \chi
\cdot \nabla \varphi_p + \sum_{0 \leq i,j\leq p\atop{i+j\geq p}} \nabla \varphi_i \cdot \nabla
\varphi_j \Big \}$$   
$$+ t^{-\gamma} \Big \{ g_0(\widetilde{q} , \widetilde{q}) + 2 g_0(\widetilde{q}, V) + g_0 (V -
W_p, V + W_p) + \sum_{0\leq i, j\leq p \atop{i+j\geq p+1}} g_0(w_i \ w_j) \Big \} \eqno(6.21)$$

\noi where the last bracket is obtained by rewriting
$$g_0(w_{t_0}, w_{t_0}) - \sum_{i+j\leq p} g_0(w_i , w_j) \quad .$$

\noi We estimate $(\widetilde{q}, \widetilde{\psi})$ by Lemmas 3.2 and 3.3, especially (3.6)
(3.12) for $\widetilde{q}$ and (3.7) (3.8) (3.13) for $\widetilde{\psi}$ and we obtain 
$$\Big | \partial_t |\widetilde{q} |_k \Big | \leq C \ t^{-2} \left \{ |V|_{k+2} +
|\varphi_{t_0}|_{\ell} \ |\widetilde{q}|_k + |\widetilde{\psi} + \varphi_p + \chi |_{\ell} \
|V|_{k+1} \right \} \eqno(6.22)$$  
$$\Big | \partial_t |\widetilde{\psi} |_{\ell} \Big | \leq C \ t^{-2} \Big \{
|\widetilde{\psi}|_{\ell}^2 + |\widetilde{\psi}|_{\ell} \ |\phi_p + \chi |_{\ell + 1} +
|\chi|_{\ell + 1}^2  + |\varphi_p|_{\ell +1} \ |\chi|_{\ell +1} + \sum_{0\leq i, j \leq p
\atop{i+j\geq p}} |\varphi_i|_{\ell + 1} \ |\varphi_j|_{\ell + 1} \Big \}$$ $$+ C \
t^{-\gamma} \Big \{ |\widetilde{q}|_k^2 + |\widetilde{q}|_k \ |V|_k + |V - W_p|_k \ |V + W_p|_k +
\sum_{0 \leq i,j \leq p\atop{i+j \geq p+1}} |w_i|_k \ |w_j|_k \Big \} \ . \eqno(6.23)$$  

\noi From Propositions 5.1 and 5.4, it follows that there exist $a$ and $b$ depending on $(a_+,
b_+)$, such that the following estimates hold.
\begin{eqnarray*}
&|V|_{k+2} \leq a \quad , \quad |W_p|_k \leq a \quad , \quad |V - W_p|_k \leq a \ Q_p \quad , &\\
&\hskip - 3.3 truecm |w_j|_k \leq a \ Q_{j-1} &\qquad \hbox{for} \ 1 \leq j \leq p \quad , \\
&|\phi_p|_{\ell + 1} \leq b \ h_0 \quad , \quad |\varphi_j|_{\ell + 1} \leq b \ N_j &\qquad
\hbox{for} \ 0 \leq j \leq p \quad , \\
&|\chi |_{\ell + 1} \leq b \leq C \ b \ N_p \leq C \ b \ \gamma^{-p} \ h_0 &\qquad \hbox{for} \  t
\geq 2 \quad .\end{eqnarray*} 

\noi We now define $y = |\widetilde{q}|_k$ and $z = |\widetilde{\psi}|_{\ell}$. Using the previous
estimates, we obtain from (6.22) (6.23)
$$|\partial_t y | \leq C \ t^{-2} \Big \{ a + (z + b \ h_0)y + (z + b \ N_p)a \Big \} 
\eqno(6.24)$$
$$|\partial_t z | \leq C \ t^{-2} \Big \{ (z + b \ h_0)z + b^2 \sum_{0 \leq i,j \leq p\atop{i+j
\geq p}} N_i \ N_j \Big \} + C \ t^{-\gamma} \Big \{ y(y + a) + a^2 \ Q_p + a^2 \sum_{0 \leq i,j
\leq p-1\atop{i+j \geq p-1}} Q_i \ Q_j \Big \} \ . \eqno(6.25)$$

\noi In the last bracket in (6.25), the terms in $a^2$ are absent for $p = 0$ since in that case
$V = W_0 = w_0 = w_+$. \par

We next estimate $y$ and $z$ from (6.24) (6.25), taking $C = 1$ for the rest of the proof. We
distinguish two cases. \\

\noi {\bf Case t $\geq$ t$_{\bf 0}$.} Let $\bar{t} > t_0$ and define $Y = Y(\bar{t}) = \ \parallel y;
L^{\infty} ([t_0, \bar{t}])\parallel$ and $Z = Z(\bar{t}) =$\break \noindent $\parallel h_0^{-1} z;
L^{\infty} ([t_0, \bar{t}])\parallel$. Then for all $t \in [t_0, \bar{t}\, ]$
$$\partial_t \ y \leq t^{-2} \Big \{ a + (Z + b) Y h_0 + a\ Z \ h_0 + a\ b \ N_p \Big \}
\eqno(6.26)$$
$$\partial_t \ z \leq t^{-2} \Big \{ (Z + b) Z \ h_0^2 + b^2 \sum_{0 \leq i,j\leq p\atop{i+j \geq
p}} N_i \ N_j \Big \} + t^{-\gamma} \Big \{ Y (Y + a) + a^2 \ Q_p + a^2 \sum_{0 \leq i,j \leq
p-1\atop{i+j \geq p-1}} Q_i \ Q_j \Big \} \ . \eqno(6.27)$$

Integrating (6.26) between $t_0$ and $t$ and using (3.35) we obtain 
$$y \leq a \ t_0^{-1} + (Z + b) Y \ h(t_0) + a\ Z\ h(t_0) + a\ b \ Q_p (t_0) \eqno(6.28)$$

\noi and therefore
$$Y \leq (Z + b) Y \ h(t_0) + a \ Z \ h(t_0) + a\ B_1(t_0) \eqno(6.29)$$

\noi where
$$B_1(t_0) = t_0^{-1} + b \ Q_p(t_0) \leq \left ( Q_p (1)^{-1} + b \right ) Q_p (t_0) \quad .
\eqno(6.30)$$

Integrating (6.27) between $t_0$ and $t$, we obtain similarly
$$z \leq (Z + b) Z \ h_0(t) \ h(t_0) + Y(Y + a) \ h_0(t) + C(b^2 + a^2) \ h_0(t) \ Q_p(t_0) \quad
. \eqno(6.31)$$

\noi Here we have used the following relations
$$\int_{t_0}^t dt_1 \ t_1^{-2} \ h_0^2(t_1) \leq h_0(t) (h(t_0) - h(t))$$
$$\int_{t_0}^t dt_1 \ t_1^{-2} \ N_i(t_1) \ N_j(t_1) \leq \int_{t_0}^t dt_1 \ t_1^{-2} \ h_0 (t_1)
\ N_{i+j}(t_1)$$
$$\leq h_0(t) \int_{t_0}^t dt_1 \ t_1^{-2} \ N_{i+j}(t_1) = h_0(t) \left ( Q_{i+j} (t_0) -
Q_{i+j}(t) \right )$$

\noi by (3.38) (3.35),
$$\int_{t_0}^t dt_1 \ t_1^{-\gamma} \ Q_p (t_1) \leq Q_p(t_0) \left ( h_0(t) - h_0(t_0) \right )$$

\noi by (3.45), and
$$\int_{t_0}^t dt_1 \ t_1^{-\gamma} \ Q_i(t_1) \ Q_j(t_1) \leq \int_{t_0}^t dt_1 \ t_1^{-\gamma} \
h(t_1)\ Q_{i+j}(t_1) \leq 2 Q_{i+j+1}(t_0) \left ( h_0(t) - h_0(t_0) \right )$$

\noi by (3.40) and (3.46). \par

>From (6.31) we obtain
$$Z \leq (Z + b) Z \ h(t_0) + Y(Y+a) + B_2(t_0) \eqno(6.32)$$

\noi with
$$B_2 = C ( b^2 + a^2 ) \ Q_p(t_0) \quad . \eqno(6.33)$$

Now (6.29) (6.32) define a closed subset ${\cal R}$ of ${I \hskip - 1 truemm R}^+ \times {I \hskip - 1 truemm
R}^+$ in the $(Y, Z)$ variables, containing the point $(0,0)$, and $(Y,Z)$ is a continuous
function of $\bar{t}$ starting from that point for $\bar{t} = t_0$. If we can find an open region
${\cal R}_1$ of ${I \hskip - 1 truemm R}^+ \times {I \hskip - 1 truemm R}^+$ containing $(0,0)$
and such that $\overline{{\cal R} \cap {\cal R}_1} \subset {\cal R}_1$, then $(Y, Z)$ will remain
in ${\cal R} \cap {\cal R}_1$ for all time, because $\overline{{\cal R} \cap {\cal R}_1}$ is both
open and closed in ${\cal R}$. We first take $t_0$ sufficiently large so that 
$$4b \ h(t_0) \leq 1 \quad , \quad 16a^2 \ h(t_0) \leq 1 \quad , \quad 4B_1(t_0) \leq 1
\eqno(6.34)$$

\noi and we choose for ${\cal R}_1$ the region $4Z \ h(t_0) < 1$. From (6.29) (6.32) (6.34) it
follows that in $\overline{{\cal R} \cap {\cal R}_1}$

$$\left \{ \begin{array}{l} Y \leq 2a \left ( Z \ h(t_0) + B_1(t_0) \right ) \leq a \\
\\ Z \leq 4a \ Y + 2 B_2(t_0) \end{array} \right .$$

\noi and therefore

$$\left \{ \begin{array}{l} Y \leq 4a \ B_1(t_0) + 2B_2(t_0) \ h(t_0)  \\
\\ Z \leq 4 \left ( 4a^2 \ B_1(t_0) +  B_2(t_0) \right ) \end{array} \right . \eqno(6.35)$$

\noi so that the condition $4Zh(t_0) < 1$ is implied by 
$$16 \left ( 4a^2 \ B_1(t_0) + B_2(t_0) \right ) \ h(t_0) < 1 \quad . \eqno(6.36)$$

 The estimates (6.15) (6.16) with $V$ and $\chi$ now follow from (6.35), while the conditions
(6.34) (6.36) reduce to the form (6.13). \par

The estimates (6.15) (6.16) with $W_p$ and $\psi_+$ follow from the previous ones, from (5.43)
(5.44) and from the fact that
$$Q_p(t) \ h_0(t) \ h(t)^{-1} = \left ( t \ Q_p(t) \right ) \left ( h_0(t) \ t^{-1} \ h(t)^{-1}
\right )$$

\noi is an increasing function of $t$, so that
$$h(t) \leq h(t_1) \ Q_p(t_1)^{-1} \ h_0(t_1)^{-1} \ Q_p(t_0) \ h_0(t)$$

\noi for any (fixed) $t_1 \leq t_0$. \\

\noi {\bf Case t $\leq$ t$_{\bf 0}$.} Let $\bar{t} < t_0$ and define $Y = Y(\bar{t}) = \parallel
Q_p^{-1} \ y; L^{\infty} ([\bar{t}, t_0 ])\parallel$ and $Z = Z(\bar{t}) =$\break \noindent
$\parallel P_p^{-1} \ z; L^{\infty} ([\bar{t}, t_0 ])\parallel$. It then follows from (6.24) (6.25)
that for all $t \in [\bar{t},t_0 ]$

$$|\partial_t \ y| \leq t^{-2} \Big \{ a + Z\ Y \ P_p \ Q_p + b\ Y \ h_0 \ Q_p + a\ Z \ P_p +
a\ b \ N_p\Big \} \eqno(6.37)$$
$$|\partial_t \ z| \leq t^{-2} \Big \{ Z^2 \ P_p^2 + b \  Z \ h_0 \ P_p + b^2 \sum_{0 \leq i,j\leq
p\atop{i+j \geq p}} N_i \ N_j \Big \}$$ $$+ t^{-\gamma} \Big \{ Y^2 \ Q_p^2 + a \ Y \  Q_p + a^2
\ Q_p + a^2 \sum_{0 \leq i,j \leq p-1\atop{i+j \geq p-1}} Q_i \ Q_j \Big \} \ . \eqno(6.38)$$

Integrating (6.37) between $t$ and $t_0$, using (3.32) (3.35) and 
\begin{eqnarray*}
&&\int_t^{t_0} dt_1 \ t_1^{-2} \ P_p (t_1) \ Q_p(t_1) \leq Q_p (t) \left ( R_p(t) -
R_p(t_0) \right ) \\
&& \int_t^{t_0} dt_1 \ t_1^{-2} \ h_0 (t_1) \ Q_p(t_1) \leq Q_p (t) \left ( h(t) - h(t_0) \right
)\end{eqnarray*}

\noi we obtain
$$y \leq a \ t^{-1} + Z \ Y \ R_p(t) \ Q_p(t)  + b \ Y \ h(t) \ Q_p (t) + a\ Z \ R_p(t) + a
\ b \ Q_p(t) \eqno(6.39)$$

\noi and therefore by (3.47)
$$ Y \leq a \ c + Z\ Y \ R_p(t) + b \ Y \ h(t) + C_p \ a \ Z \ h(t) + a \ b \eqno(6.40)$$

\noi where $c = Q_p(1)^{-1}$. \par

We integrate similarly (6.38) between $t$ and $t_0$. We use the relations
\begin{eqnarray*}
\int_t^{t_0} dt_1 \ t_1^{-2} \ P_p^2(t_1) &\leq& P_p(t) \ \left ( R_p(t) - R_p(t_0)\right ) \leq P_p
(t) \ R_p(t) \\
\int_t^{t_0} dt_1 \ t_1^{-2} \ h_0(t_1) \ P_p(t_1) &\leq& P_p (t) \left ( h(t) - h(t_0)\right ) \\
\int_t^{t_0} dt_1 \ t_1^{-2} \ N_i(t_1) \ N_j(t_1) &\leq& \int_t^{t_0} dt_1 \ t_1^{-2} \ h_0(t_1) \
N_{i+j}(t_1) \\
&\leq& P_{i+j} (t) - P_{i+j}(t_0)
 \end{eqnarray*}

\noi by (3.38) (3.37),
$$\int_t^{t_0} dt_1 \ t_1^{-\gamma} \ Q_p^2(t_1) \leq Q_p(t) \int_t^{t_0} dt_1 \ t_1^{-\gamma} \
Q_p(t_1) \leq Q_p(t) \ P_p(t)$$

\noi by (3.42) and
$$\int_t^{t_0} dt_1 \ t_1^{-\gamma} \ Q_i(t_1) \ Q_j(t_1) \leq P_{i+j+1}(t)$$

\noi by (3.40) (3.41) (3.42). We obtain
$$z \leq Z^2 \ P_p(t) \ R_p(t) + b\ Z \ P_p(t) \ h(t) + Y^2 \ P_p(t) \ Q_p(t) + a \ Y \ P_p(t) + C
( a^2 + b^2 ) P_p(t) \eqno(6.41)$$

\noi and therefore
$$ Z \leq Z^2 \ R_p(t) + b \ Z \ h(t) + Y^2 \ Q_p(t) + a\ Y + C  ( a^2 + b^2   ) \quad .
\eqno(6.42)$$

\noi We now take $t$ sufficiently large so that $bh(t) \leq 1/4$ and we proceed as in the case $t
\geq t_0$ by taking for ${\cal R}_1$ the strip defined by $ZR_p(t) < 1/4$, $C_pZh(t) < b$,
thereby obtaining from (6.40) (6.42) 
$$\left \{ \begin{array}{l} Y \leq 2a (2b + c) \\ \\ Z \leq 2C \left ( a^2 + b^2 \right ) + 4a^2
(2b + c) \left ( 1 + 2 (2b + c) Q_p(1) \right ) \quad . \end{array} \right . \eqno(6.43)$$

\vskip 3 truemm
\noi The conditions $bh \leq 1/4$, $ZR < 1/4$, $Zh < b$ are then satisfied for $t \geq T$ with $T$
defined by a condition of the form (6.14), where the singular factor $((p+2)\gamma - 1)^{-1}$
comes from $C_p$. The estimates (6.17) (6.18) with $V$ and $\chi$ follow from (6.43), and the
analogous estimates with $W_p$ and $\psi_+$ follow therefrom and from (5.43) (5.44). \par

Finally the estimates (6.19) follow from (6.15) (6.16) (6.17) (6.18) and from Proposition 5.1.
\par \nobreak
\hfill $\sq$

We can now take the limit $t_0 \to \infty$ of the solution $(w_{t_0}, \varphi_{t_0})$ constructed
in Proposition 6.2, for fixed $(w_+, \psi_+)$. \\

\noi {\bf Proposition 6.3.} {\it Let $(k, \ell )$ be an admissible pair and let $p$ be an integer
such that $(p + 2)\gamma > 1$. Let $w_+ \in H^{k+(p+1)\vee 2}$ and $\psi_+ \in Y^{\ell + 1}$.
Let $\phi = \phi_{p-1}$ be defined by (5.11) and Proposition 5.1, so that $h_0^{-1}\phi \in
({\cal C} \cap L^{\infty})([1, \infty ), Y^{\ell+2})$. Let $V$ and $\chi$ be the solutions of
(5.12) and (5.39) respectively, obtained in Proposition 5.4 so that $(V, \chi ) \in ({\cal C}
\cap L^{\infty})([1, \infty ), H^{k+2} \oplus Y^{\ell + 1})$ and let $a_+$, $b_+$ be defined by
(6.12). Then}\par
{\it (1) There exists $T$, $1 \leq T < \infty$, depending only on $(\gamma , p, a_+, b_+)$ and
there exists a unique solution $(w, \varphi )$ of the system (4.1) (4.2) in the interval $[T,
\infty )$ such that $(w, h_0^{-1} \varphi ) \in ({\cal C} \cap L^{\infty})([T, \infty ), H^{k}
\oplus Y^{\ell})$ and such that the following estimates hold for all $t \geq T$.}
$$|w(t) - V(t)|_k \vee \Big |w(t) - W_p(t) \Big |_k \leq  A\left ( a_+, b_+ \right ) \ Q_p(t)
\eqno(6.44)$$
$$\Big |\varphi (t) - \phi_p(t) - \chi (t) \Big |_{\ell} \vee \Big |\varphi (t) - \phi_p(t) - \psi_+
\Big |_{\ell} \leq  A\left ( a_+, b_+ \right ) \ P_p(t) \eqno(6.45)$$  
  $$|w(t)|_k \leq A\left ( a_+, b_+ \right ) \quad , \quad  |\varphi (t)|_{\ell} \leq  A\left
( a_+, b_+ \right ) \ h_0(t) \eqno(6.46)$$

\noi {\it One can define $T$ by a condition of the type (6.14).} \par

{\it (2) Let $(w_{t_0}, \varphi_{t_0})$ be the solution of the system (4.1) (4.2) constructed in
Proposition 6.2 for $t_0 \geq T_0 \vee T$ and such that $(w_{t_0}, h_0^{-1} \varphi_{t_0}) \in
({\cal C} \cap L^{\infty})([T, \infty ), H^k \oplus Y^{\ell})$. Then $(w_{t_0} , \varphi_{t_0})$
converges to $(w, \varphi )$ in norm in $L^{\infty}(J, H^{k-1} \oplus Y^{\ell - 1})$ and in the
weak-$*$ sense in $L^{\infty}(J, H^k \oplus Y^{\ell})$ for any compact $J \subset [T, \infty )$,
and in the weak-$*$ sense in $H^k \oplus Y^{\ell}$ pointwise in $t$.} \par
{\it (3) The map $(w_+, \psi_+) \to (w, \varphi )$ defined in Part (1) is continuous on the
bounded sets of $H^{k+(p+1)\vee 2} \oplus Y^{\ell + 1}$ from the norm topology of $(w_+, \psi_+)$
in $H^{k+p-1} \oplus Y^{\ell - 1}$ to the norm topology of $(w, \varphi )$ in $L^{\infty}(J,
H^{k-1} \oplus Y^{\ell - 1})$ and to the weak-$*$ topology in $L^{\infty}(J, H^k \oplus Y^{\ell})$
for any compact interval $J \subset [T, \infty )$, and to the weak-$*$ topology in $H^k \oplus
Y^{\ell}$ pointwise in $t$.} \\

\noi {\bf Proof.} {\bf Parts (1) and (2)} will follow from the convergence of $(w_{t_0},
\varphi_{t_0})$ when $t_0 \to \infty$ in the topologies stated in Part (2). Let $T_0 \vee T \leq
t_0 \leq t_1$. From (6.17) (6.18) it follows that 
$$\Big | w_{t_1}(t_0) - w_{t_0}(t_0) \Big |_k = \Big | w_{t_1}(t_0) - V(t_0) \Big |_k \leq A \
Q_p(t_0) \eqno(6.47)$$
$$\Big | \varphi_{t_1}(t_0) - \varphi_{t_0}(t_0) \Big |_{\ell} = \Big | \varphi_{t_1}(t_0) -
\phi_p(t_0) - \chi(t_0) \Big |_{\ell} \leq A \ P_p(t_0) \eqno(6.48)$$

\noi We now estimate $(w_{t_0} - w_{t_1}, \varphi_{t_0} - \varphi_{t_1})$ in $H^{k-1} \oplus
Y^{\ell - 1}$ for $t \leq t_0$. Let
$$y = \Big | w_{t_0} - w_{t_1} \Big |_{k-1} \quad , \quad z = \Big | \varphi_{t_0} - \varphi_{t_1}
\Big |_{\ell - 1} \quad . \eqno(6.49)$$

\noi From (6.19) and Lemma 3.3, it follows that $y$ and $z$ satisfy the system (4.10). Integrating
that system for $t \leq t_0$ with initial data at $t_0$, we obtain from Lemma 4.1 

$$\left \{ \begin{array}{l} y(t) \leq A \left ( y(t_0) + t^{-1} \left ( z(t_0) + y (t_0) h_0(t_0)
\right ) \right ) \\ \\ z(t) \leq A \left ( z(t_0) + y(t_0) h_0(t_0) \right ) \quad .  
\end{array}\right . \eqno(6.50)$$
\vskip 3 truemm

\noi From (6.47) (6.48) (6.50) and from the fact that $P_p(t_0)$ and $Q_p(t_0) h_0(t_0)$ tend to
zero when $t_0 \to \infty$, it follows that there exists $(w, \varphi ) \in {\cal C}([T, \infty ),
H^{k-1} \oplus Y^{\ell - 1})$ such that $(w_{t_0}, \varphi_{t_0})$ converges to $(w, \varphi )$ in
$L^{\infty}(J, H^{k-1} \oplus Y^{\ell - 1})$ for all compact intervals $J \subset [T, \infty )$.
>From that convergence, from (6.17) (6.18) (6.19) and standard compactness arguments, it follows
that $(w, h_0^{-1}\varphi ) \in ({\cal C}_{w*} \cap L^{\infty})([T, \infty ), H^k \oplus
Y^{\ell})$, that $(w, \varphi )$ satisfies the estimates (6.44) (6.45) (6.46) for all $t \geq T$,
and that $(w_{t_0}, \varphi_{t_0})$ converges to $(w, \varphi )$ in the other topologies
considered in Part (2). Furthermore, $(w, \varphi )$ satisfies the system (4.1) (4.2) and by
Proposition 4.1, part (1), $(w, \varphi ) \in {\cal C}([T, \infty ), H^k \oplus Y^{\ell})$.
Finally, uniqueness of $(w, \varphi )$ under the conditions (6.44) (6.45) follows from
Proposition 4.3 and from the fact that $P_p(t)$ and $Q_p(t) h_0(t)$ tend to zero when $t \to
\infty$.  \\ 
  
\noi {\bf Part (3).} Let $(w_+, \psi_+)$ and $(w'_+ , \psi '_+)$ belong to a fixed bounded set
of $H^{k+(p+1)\vee 2} \oplus Y^{\ell + 1}$. Let $(W_p, \phi_p)$ and $(W'_p, \phi'_p)$ be
the associated functions defined by (5.11) and Proposition 5.1 and let $(w, \varphi )$ and
$(w', \varphi ')$ be the associated solutions of the system (4.1) (4.2) defined in Part (1). We
assume that $(w'_+, \psi '_+)$ is close to $(w_+, \psi_+)$ in the sense that
$$|w_+ - w'_+|_{k+p-1} \leq \varepsilon \eqno(6.51)$$
$$|\psi_+ - \psi '_+ |_{\ell - 1} \leq \varepsilon_0 \quad . \eqno(6.52)$$

\noi We now take $t_0 > T$ and we estimate $(w - w', \varphi - \varphi ')$ in $H^{k-1} \oplus
Y^{\ell - 1}$ for $t \leq t_0$. Let
$$y = |w - w'|_{k-1} \qquad , \qquad z = |\varphi - \varphi '|_{\ell - 1} \quad . \eqno(6.53)$$

\noi From (4.1) (4.2) and Lemma 3.3, it follows that $(y, z)$ satisfy the system (4.10).
Integrating that system between $t_0$ and $t$ yields the estimate (6.50) for $(y, z)$ defined
by (6.53). From (6.44) (6.45) we obtain

$$\left \{ \begin{array}{l} y(t_0) \leq A \ Q_p(t_0) + \Big | W_p(t_0) - W'_p(t_0) \Big |_{k-1}  \\
\\ z(t_0) \leq A \ P_p(t_0) + \varepsilon_0 + \Big | \phi_p(t_0) - \phi '_p (t_0) \Big |_{\ell - 1}
\quad .   \end{array}\right . \eqno(6.54)$$

\noi From estimates similar to (5.8) (5.9), and from (6.51) (6.52) it then follows that
$$\left \{ \begin{array}{l} y(t_0) \leq A \left ( Q_p(t_0) + \varepsilon \right )
\\ \\ z(t_0) \leq A \left ( P_p(t_0) + \varepsilon \ h_0(t_0) \right ) + \varepsilon_0 \quad .
\end{array}\right . \eqno(6.55)$$

\noi We now choose $t_0$ so that $Q_p(t_0) = \varepsilon$. Substituting (6.55) with that choice
into (6.50) and using the asymptotic behaviour of $Q_p$ and $P_p$ for large $t$, we obtain
 $$\left \{ \begin{array}{l} y(t) \leq A \left ( m(\varepsilon ) + t^{-1} \left ( \varepsilon_0 +
m(\varepsilon ) \right ) \right ) \\ \\ z(t) \leq A \left ( \varepsilon_0 + m(\varepsilon ) \right )
\end{array}\right . \eqno(6.56)$$

\noi where
$$m (\varepsilon ) = \left \{ \begin{array}{ll} \varepsilon ^{((p+2)\gamma -
1)/(p+1)\gamma} &\qquad \hbox{for} \ (p+1)\gamma < 1 \\
\\
\varepsilon^{\gamma} \ {\rm Log} \ \varepsilon &\qquad \hbox{for} (p+1) \gamma = 1
\\ \\ \varepsilon^{\gamma} &\qquad \hbox{for} (p+1)\gamma > 1 \quad . \end{array}\right .
\eqno(6.57)$$

This implies the (uniform H\"older) continuity of $(w, \varphi )$ as a function of $(w_+, \psi_+)$
in the norm topology of $L^{\infty}(J,H^{k-1} \oplus Y^{\ell - 1})$ for all compact intervals $J
\subset [T, \infty )$. The other continuities follow therefrom and from the boundedness of $(w,
h_0^{-1} \varphi )$ in $L^{\infty}([T, \infty ), H^k \oplus Y^{\ell})$ by standard compactness
arguments. \par \nobreak
\hfill $\sq$ 

\section{Asymptotics and wave operators for $u$}
\hspace*{\parindent} In this section we complete the construction of the wave operators for the equation (1.1) and we
derive asymptotic properties of solutions in their range. The construction relies in an essential
way on those of Section 6, esp. Proposition 6.3, and will require a discussion of the gauge
invariance of those constructions. \par

We first define the wave operator for the auxiliary system (4.1) (4.2).\\

\noi {\bf Definition 7.1.} We define the wave operator $\Omega_0$ as the map
$$\Omega_0 : (w_+, \psi_{+} ) \to (w, \varphi ) \eqno(7.1)$$

\noi from $H^{k+(p+1)\vee 2} \oplus Y^{\ell + 1}$ to the space of $(w, \varphi )$ such that $(w,
h_0^{-1} \varphi ) \in ({\cal C} \cap L^{\infty}) ([T, \infty ), H^k \oplus Y^{\ell})$
for some $T$, $1 \leq T < \infty$, where $(w, \varphi )$ is the solution of the system (4.1)
(4.2) obtained in Proposition 6.3, part (1). \\

Before defining the wave operators for $u$, we now study the gauge invariance of $\Omega_0$, which
plays an important role in justifying that definition, as was explained in Section 2. For
that purpose we need some information on the Cauchy problem for the equation (1.1) at
finite times. In addition to the operators $M = M(t)$ and $D = D(t)$ defined by (2.4)
(2.5), we introduce the operator
$$J = J(t) = x + it \nabla \quad , \eqno(7.2)$$
\noi the generator of Galilei transformations. The operators $M$, $D$, $J$ satisfy the
commutation relation
$$i \ M \ D \ \nabla = J \ M \ D \quad . \eqno(7.3)$$
\noi For any interval $I \subset [1, \infty )$ and any nonnegative integer $k$, we define the
space
$$\begin{array}{ll} {\cal X}^k(I) &= \left \{ u : D^* M^* u \in
{\cal C} (I, H^k) \right \}  \\
& \\
&= \left \{ u : <J(t)>^k u \in {\cal C} (I, L^2) \right \} 
\end{array} \eqno(7.4)$$ 
\noi where $<\lambda > = (1 + \lambda^2)^{1/2}$ for any real number or self-adjoint
operator $\lambda$ and where the second equality follows from (7.3). Then \cite{3r} \\

\noi {\bf Proposition 7.1.} {\it Let $k$ be a positive integer and let $0 < \mu < n \wedge 2k$.
Then the Cauchy problem for the equation (1.1) with initial data $u(t_0) = u_0$ such that
$<J(t_0)>^k$ $u_0 \in L^2$ at some initial time $t_0 \geq 1$ is locally well posed in
${\cal X}^k(\cdot )$, namely} \par
{\it (1) There exists $T > 0$ such that (1.1) has a unique solution with initial data $u(t_0) = u_0$
in ${\cal X}^k ([1 \vee (t_0 - T), t_0 + T])$.} \par
{\it (2) For any interval $I$, $t_0 \in I \subset [1 , \infty )$, (1.1) with initial data $u(t_0) =
u_0$ has at most one solution in ${\cal X}^k(I)$.} \par
{\it (3) The solution of Part (1) depends continuously on $u_0$ in the norms considered there.} \\

We come back from the system (4.1) (4.2) to the equation (1.1) by reconstructing $u$ from $(w,
\varphi )$ by (2.7) and accordingly we define the map

$$\Lambda : (w, \varphi ) \to u = M \ D \exp (-i \varphi ) w \quad . \eqno(7.5)$$

\noi It follows immediately from Lemma 3.1 that the map $\Lambda$ satisfies the following
property.\\

\noi {\bf Lemma 7.1.} {\it The map $\Lambda$ defined by (7.5) is bounded and continuous from ${\cal
C} (I, H^k \oplus Y^{\ell})$ to ${\cal X}^k(I)$ for any admissible pair $(k, \ell )$ and any
interval $I \subset [1, \infty )$.} \\

We now give the following definition. \\

\noi {\bf Definition 7.2.} Let $(k, \ell )$ be an admissible pair and let $(w, \varphi )$ and
$(w', \varphi ')$ be two solutions of the system (4.1) (4.2) in ${\cal C} (I, H^k \oplus
Y^{\ell})$ for some interval $I \subset [1, \infty)$. We say that $(w, \varphi )$ and $(w',
\varphi ')$ are gauge equivalent if they give rise to the same $u$, namely if $\Lambda (w, \varphi
) = \Lambda (w', \varphi ')$, or equivalently if
$$\exp ( - i \varphi (t) ) \ w(t) = \exp (- i \varphi '(t)) \ w'(t)
\eqno(7.6)$$ 

\noi for all $t \in I$. \par

A sufficient condition for gauge equivalence is given by the following Lemma. \\   

\noi {\bf Lemma 7.2.} {\it Let $(k, \ell )$ be an admissible pair and let $(w, \varphi )$
and $(w', \varphi ')$ be two solutions of the system (4.1) (4.2) in ${\cal C} (I, H^k
\oplus Y^{\ell})$. In order that $(w, \varphi )$ and $(w', \varphi ')$ be gauge equivalent,
it is sufficient that (7.6) holds for one $t \in I$.} \\

\noi {\bf Proof.} An immediate consequence of Lemma 7.1, of Proposition 7.1, part (2), and
of the fact that $(k, \ell)$ admissible implies $k > 1 + \mu /2$. \par \nobreak
\hfill $\sq$ \par

The gauge covariance properties of $\Omega_0$ will be expressed by the following two propositions.
\\

\noi {\bf Proposition 7.2.} {\it Let $(k, \ell )$ be an admissible
pair. Let $(w, \varphi )$ and $(w', \varphi ')$ be two solutions of the system (4.1) (4.2)
such that $(w, h_0^{-1}\varphi), (w', h_0^{-1}\varphi ') \in ({\cal C} \cap
L^{\infty})([T, \infty ), H^k \oplus Y^{\ell})$ for some $T \geq 1$, and assume that $(w,
\varphi )$ and $(w', \varphi ')$ are gauge equivalent. Then} \par

{\it (1) There exists $\sigma \in Y^{\ell - 1}$ such that $\varphi ' (t) - \varphi (t)$
converges to $\sigma$ when $t \to \infty$ strongly in $Y^{\ell - 2}$ and in the weak-$*$ sense in
$Y^{\ell - 1}$. The following estimates hold~:}
$$\parallel \varphi ' (t) - \varphi (t) - \sigma ; Y^{\ell - 2} \parallel \ \leq A \ h(t)
\eqno(7.7)$$ \noi {\it for some constant $A$ depending on $T$ and on the norms of
$h_0^{-1} \varphi$, $h_0^{-1} \varphi '$ in $L^{\infty} (\cdot , Y^{\ell})$, with
the exception of the case $n$ even, $\ell = n/2 + 1$ where the $L^{\infty}$ norm of
$\nabla \varphi$ satisfies only}
$$\parallel \nabla \varphi '(t) - \nabla \varphi (t) - \nabla \sigma \parallel_{\infty}
\ \leq A \ h(t)^{1/2} \quad . \eqno(7.8)$$

{\it (2) Let $w_+$ and $w'_+$ be the limits of $w(t)$ and $w'(t)$ as $t \to \infty$, obtained
in Proposition 4.4. Then $w'_+ = w_+ \exp (- i \sigma )$.} \par

{\it (3) Let $p \geq 0$ be an integer. Assume in addition that $w_+$, $w'_+ \in H^{k+p}$ and
let $\phi_p$, $\phi '_p$ be associated with $w_+$, $w'_+$ according to (5.11) and Proposition
5.1. Assume that the following limits exist}
$$\lim_{t \to \infty} \Big ( \varphi (t) - \phi_p(t) \Big ) = \psi_+ \qquad , \qquad \lim_{t \to
\infty} \Big ( \varphi '(t) - \phi '_p(t) \Big ) = \psi '_+ \eqno(7.9)$$

\noi {\it as strong limits in $L^{\infty}$. Then $\psi '_+ = \psi_+ + \sigma$.} \\     

\noi {\bf Proof.} {\bf Part (1)} is essentially identical with Proposition I.7.2, part (1). \\

\noi {\bf Part (2).} We define $\varphi_-(t) = \varphi '(t) - \varphi (t)$ and we estimate 
$$\parallel w'_+ - w_+ \ e^{i\sigma} \parallel_2 \ \leq \ \parallel w'_+ - w'(t) \parallel_2 \ + \
\parallel w'(t) - \exp (i \varphi_- (t)) w(t) \parallel_2$$
$$+ \parallel  ( \exp (i \varphi_-(t)  ) - \exp (i \sigma )) w(t) \parallel_2 \ + \ \parallel
\exp (i \sigma ) (w(t) - w_+)\parallel_2$$
$$\leq \ \parallel w'_+ - w'(t) \parallel_2 \ + \ \parallel w(t) - w_+\parallel_2 \ + \ \parallel
\varphi_-(t) - \sigma \parallel_{\infty} \ \parallel w(t) \parallel_2 \eqno(7.10)$$

\noi by gauge invariance. The last member of (7.10) tends to zero as $t \to \infty$. \\

\noi {\bf Part (3).} By gauge invariance, namely Proposition 5.1 part (2) and Part (2) of this
proposition, $\phi '_p = \phi_p$ and therefore
$$\psi '_+ - \psi _+ = \lim_{t \to \infty} \Big ( \varphi '(t) - \varphi (t) \Big ) = \sigma \quad .
\eqno(7.11)$$
 \hfill $\sq$ \par

\noi {\bf Remark 7.1.} The additional assumptions of Proposition 7.2, part (3) are satisfied
either if $(w, \varphi )$, $(w', \varphi ')$ satisfy the assumptions of Proposition 6.1, or if $(w,
\varphi )$, $(w', \varphi ') \in {\cal R}(\Omega_0)$. We shall not consider the former case any
further. In the latter case, it follows from (7.11) that actually $\sigma \in Y^{\ell + 1}$.\\

Proposition 7.2 prompts us to make the following definition of gauge equivalence for asymptotic
states. \\

\noi {\bf Definition 7.3.} Two pairs $(w_+ , \psi_+)$ and $(w'_+, \psi '_+)$ are gauge equivalent
if $w_+ \exp (- i \psi_+) = w'_+ \exp (- i \psi '_+)$.  \\

With this definition, Proposition 7.2 implies that two gauge equivalent solutions of the system
(4.1) (4.2) in ${\cal R} (\Omega_0)$ are images of two gauge equivalent pairs of asymptotic
states. The next proposition shows that conversely two gauge equivalent pairs of asymptotic
states have gauge equivalent images under $\Omega_0$. \\

\noi {\bf Proposition 7.3.} {\it Let $(k, \ell )$ be an admissible pair and let $p$ be an integer
such that $(p+2)\gamma > 1$. Let $(w_+, \psi_+)$, $(w'_+, \psi '_+) \in H^{k+(p+1)\vee2} \oplus
Y^{\ell + 1}$ be gauge equivalent, and let $(w, \varphi )$, $(w', \varphi ')$ be their images under
$\Omega_0$. Then $(w, \varphi )$ and $(w', \varphi ')$ are gauge equivalent.} \\

\noi {\bf Proof.} Let $t_0$ be sufficiently large and let $(w_{t_0}, \varphi_{t_0})$ and
$(w'_{t_0}, \varphi '_{t_0})$ be the solutions of the system (4.1) (4.2) constructed by
Proposition 6.2. From the initial conditions
$$w_{t_0}(t_0) = V(t_0) \qquad , \qquad w'_{t_0} (t_0) = V'(t_0) \quad ,$$
$$\varphi_{t_0}(t_0) = \phi_p(t_0) + \chi (t_0) \quad , \quad \varphi '_{t_0}(t_0) = \phi '_p(t_0)
+ \chi ' (t_0) \quad , $$
\noi from the fact that $\phi_p = \phi '_p$ by Proposition 5.1 part (2) and that $V \exp (- i \chi
) = V' \exp (- i \chi ')$ by Proposition 5.4, part (2), it follows that
$$w_{t_0}(t_0) \exp  (- i \varphi_{t_0}(t_0) ) = w'_{t_0} (t_0) \exp ( - i
\varphi '_{t_0} (t_0) )$$

\noi and therefore by Lemma 7.2, $(w_{t_0}, \varphi_{t_0})$ and $(w'_{t_0}, \varphi '_{t_0})$ are
gauge equivalent, namely
$$w_{t_0} (t) \exp  ( - i \varphi_{t_0}(t)  ) = w'_{t_0}(t) \exp ( - i \varphi
'_{t_0}(t)  ) \eqno(7.12)$$

\noi for all $t$ for which both solutions are defined. \par

We now take the limit $t_0 \to \infty$ for fixed $t$ in (7.12). By Proposition 6.3 part (2), for
fixed $t$, $(w_{t_0}, \varphi_{t_0})$ and $(w'_{t_0}, \varphi '_{t_0})$ converge respectively to
$(w, \varphi )$ and $(w', \varphi ')$ in $H^{k-1} \oplus Y^{\ell - 1}$. By Lemma 3.1, one can take
the limit $t_0 \to \infty$ in (7.12), thereby obtaining (7.6), so that $(w, \varphi )$ and $(w',
\varphi ')$ are gauge equivalent. \par \nobreak \hfill $\sq$ \par

We can now define the wave operators for $u$. We recall from the heuristic discussion of Section 2
that we want to exploit the operator $\Omega_0$ defined in Definition 7.1, reconstruct $u$ through
the map $\Lambda$ defined by (7.5) and eliminate the arbitrariness in $\psi_+$ by fixing $\psi_+ =
0$, thereby ensuring the injectivity of the wave operator for $u$.  \\

\noi {\bf Definition 7.4.} We define the wave operator $\Omega$ as the map 
$$\Omega : u_+ \to u = \left ( \Lambda \circ \Omega_0 \right ) \left ( F \ u_+, 0 \right )
\eqno(7.13)$$

\noi from $F \ H^{k+(p+1)\vee 2}$ to ${\cal X}^k([T, \infty ))$ for some $T$, $1 \leq T < \infty$,
where $k$ is the first element of an admissible pair, and $\Omega_0$ and $\Lambda$ are defined by
Definition 7.1 and by (7.5).  \\

The fact that $\Omega$ acts between the spaces indicated follows from Proposition 6.3 and Lemma 7.1.
\\

\noi {\bf Proposition 7.4.} \par
 {\it (1) The map $\Omega$ is injective.}\par
{\it (2) If $0 \leq p \leq 2$, then ${\cal R}(\Omega ) = {\cal R} (\Lambda \circ \Omega_0)$}. \\

\noi {\bf Proof.} {\bf Part (1)} follows from the fact that $\Omega_0$ is an injective map between
gauge equivalence classes and that an equivalence class of asymptotic states contains at most one
representative with $\psi_+ = 0$.  \\

\noi {\bf Part (2)} follows from the fact that the gauge equivalence class of a given $(w_+,
\psi_+)$ actually contains an element with $\psi_+ = 0$, namely $(w_+ \exp (- i \psi_+), 0)$, by
Lemma 3.1. \par \nobreak \hfill $\sq$ \par

\noi {\bf Remark 7.2.} Part (2) of Proposition 7.4 does not extend to the case $p \geq 3$ because
in that case $(w_+, \psi_+) \in H^{k+p+1} \oplus Y^{\ell + 1}$ does not imply that $w_+ \exp (- i
\psi_+) \in H^{k+p+1}$, so that the gauge equivalence class of a given $(w_+, \psi_+)$ need not
contain an element with $\psi_+ = 0$. \\

We now collect the information obtained for the solutions of the equation (1.1) so far
cons\-truc\-ted. The main result of this paper can be stated as follows. \\

\noi {\bf Proposition 7.5.} {\it Let $n \geq 3$, $0 < \mu \leq n - 2$ and $0 < \gamma \leq 1$. Let
$(k, \ell)$ be an admissible pair. Let $p \geq 0$ be an integer with $(p+2)\gamma > 1$. Let $u_+
\in F \ H^{k+(p+1)\vee 2}$ and $a = |F \ u_+|_{k+(p+1)\vee 2}$. Let $W_p$ and $\phi_p$ be defined
by (5.11) and Proposition 5.1 with $w_+ = F \ u_+$. Then}\par
{\it (1) There exists $T$, $1 \leq T < \infty$, and there exists a unique solution $u \in {\cal
X}^k([T, \infty))$ of the equation (1.1) which can be represented as}
$$u = M \ D \exp (-i \varphi ) w$$
\noi {\it  where $(w, \varphi )$ is a solution of the system (4.1) (4.2) such that $(w,
h_0^{-1}\varphi ) \in ({\cal C} \cap L^{\infty}) ([T, \infty ), H^k \oplus Y^{\ell})$ and such that}
$$\Big | w(t) - F \ u_+ \Big |_{k-1} \ h_0(t) \to 0 \eqno(7.14)$$
$$\Big | \varphi (t) - \phi_p (t) \Big |_{\ell - 1} \to 0 \eqno(7.15)$$

\noi {\it when $t \to \infty$, where $h_0$ is defined by (3.19). The time $T$ can be defined by
(6.14) with $b_+ = 0$.} \par
{\it (2) The solution is obtained as $u = \Omega (u_+)$ where the map $\Omega$ is defined in
Definition 7.4. The map $\Omega$ is injective.} \par
{\it (3) The map $\Omega$ is continuous on the bounded sets of $F\ H^{k+(p+1)\vee 2}$ from the norm
topology in $F \ H^{k+p-1}$ for $u_+$ to the norm topology in ${\cal X}^{k-1}(I)$ and to the weak-$*$
topology in ${\cal X}^k(I)$ for $u$ for any compact interval $I \subset [T, \infty)$, and to the
weak topology in $M  D    H^k$ pointwise in $t$.} \par
{\it (4) The solution $u$ satisfies the following estimates for $t \geq T$}
$$\parallel <J(t)>^k \Big ( \exp  ( i \phi_p(t, x/t) ) u(t) - M(t) \ D(t) \ F \ u_+
\Big ) \parallel_2 \ \leq A(a) \ P_p(t) \eqno(7.16)$$

\noi {\it for some estimating function $A(a)$, where $P_p(t)$ is defined by (3.31).}\par
{\it (5) Let $r$ satisfy $0 \leq \delta (r) \equiv n/2 - n/r \leq k \wedge n/2$, $\delta (r) < n/2$
if $k = n/2$. Then $u$ satisfies the following estimate}
$$\parallel u(t) - \exp  ( - i \phi_p (t, x/t)) M(t) \ D(t) \ F \ u_+ \parallel_r \
\leq A(a) \ t^{-\delta (r)} \ P_p(t) \quad . \eqno(7.17)$$ \vskip 3 truemm

\noi {\bf Proof.} {\bf Parts (1) (2) (3)} follow from Proposition 6.3, Proposition 4.3, from
Definition 7.4, Proposition 7.4 part (1), and Lemma 7.1. \\

\noi {\bf Part (4).} From the definition (7.2) of $J(t)$, from the commutation relation (7.3) and
from Lemma 3.1, it follows that the LHS of (7.16) is estimated by 
\begin{eqnarray*}
\parallel \cdot \parallel_2 \ = \Big | \exp ( i \left ( \phi_p - \varphi
\right ) ) w - F \ u_+ \Big |_k
&\leq& |w - F \ u_+|_k + \Big | \left ( \exp ( i (\phi_p - \varphi ) ) - 1 \right )
w \Big |_k\\
&\leq& |w - F \ u_+ |_k +  | \phi_p - \varphi |_{\ell - 1} \left ( 1 + |\phi_p - \varphi
|_{\ell - 1} \right )^{k-1} \ |w|_k \quad .\end{eqnarray*}

\noi The result now follows from the estimates (4.17) and (6.45) (6.46). \\

\noi {\bf Part (5)} follows from Part (4) and from the inequality
\begin{eqnarray*}
\parallel f \parallel_r &=& t^{-\delta (r)} \parallel D^*\ M^*\ f \parallel_r \ \leq C \ t^{-\delta
(r)} \parallel < \nabla >^k \ D^* \ M^* \ f \parallel_2 \\
&=& C \ t^{-\delta (r)} \parallel <J(t)>^k \ f\parallel_2 \end{eqnarray*}  

\noi which follows from the commutation relation (7.3) and from Sobolev inequalities. \par \nobreak
\hfill $\sq$ \par
    
\noi {\bf Remark 7.3.} In (7.16) and (7.17) one could replace $MDFu_+$ by $U(t) u_+$ since $U(t)
u_+ - MDFu_+ = O(t^{-1})$ in the relevant norms. One could also replace $Fu_+$ by $W_p$, but this
would not produce any improvement in the final estimates, since the main contribution of the
dif\-fe\-ren\-ce between $u$ and its asymptotic form is that of the phase. \\

Finally, by combining Proposition 7.5 with the known results on the Cauchy problem for the equation
(1.1) at finite times, one could extend the solutions $u$ to arbitrary finite times and define more
standard wave operators $\Omega_1 : u_+ \to u(1)$ where $u = \Omega u_+$. We refer to I for the
details.

\noi {\bf Acknowledgements.} One of us (G. V.) is grateful to Professor J. C. Saut for the
hospitality at the Laboratoire d'Analyse Num\'erique et Equations aux D\'eriv\'es Partielles and to
Professor D. Schiff for the hospitality at the Laboratoire de Physique Th\'eorique. \\

\end{document}